\documentclass[11pt]{article}
\usepackage{amsmath, amssymb, amsfonts}
\usepackage{epsfig}

\def\be{\begin{equation}}
\def\ee{\end{equation}}
\def\bea{\begin{eqnarray}}
\def\eea{\end{eqnarray}}
\def\bes{\begin{eqnarray*}}
\def\ees{\end{eqnarray*}}

\def\nn{\nonumber}
\def\<{\langle}
\def\>{\rangle}
\def\lb{\label}
\def\bs{\setminus}

\def\R{{\mathbb{R}}}
\def\C{{\bf C}}
\def\Z{{\mathbb{Z}}}
\def\N{{\mathbb{N}}}

\def\Q{{\mathbb{Q}}}

\def\RP{{\mathbb{R}P^{n}}}

\def\aa{{\alpha}}
\def\bb{{\beta}}
\def\ga{{\gamma}}

\def\th{{\theta}}

\def\Lm{{\Lambda}}

\def\vf{{\varphi}}
\def\vs{{\varsigma}}

\def\ker{{\rm ker}}

\def\sgn{{\rm sgn}}

\def\hv{{\rm hv}}

\def\rank{{\rm rank}}

\def\CG{{\rm CG}}

\def\dm{{\diamond}}
\def\ol#1{\overline{#1}}  
\def\td#1{\tilde{#1}}

\def\mapright#1{\smash{\mathop{\longrightarrow}\limits^{#1}}}

\def\mapdown#1{\Big\downarrow\rlap{$\vcenter{\hbox{$\scriptstyle#1$}}$}}

\def\mapright#1{\smash{\mathop{\longrightarrow}\limits^{#1}}}

\def\mapdown#1{\Big\downarrow\rlap{$\vcenter{\hbox{$\scriptstyle#1$}}$}}

\title{Resonance identity and multiplicity of non-contractible closed geodesics on Finsler $\mathbb{R}P^{n}$}
\author{Hui Liu$^{1}$,\thanks{Partially supported by NSFC (No.11401555), Anhui Provincial Natural Science Foundation (No. 1608085QA01).
E-mail: huiliu@ustc.edu.cn.}\qquad
Yuming Xiao$^{2}$,\thanks{Supported by the Scientific Research Funds for Young Teachers of Sichuan University, Grant 2012SCU11083.
e-mail: yumingxiao@scu.edu.cn.}\\\\
$^{1}$ School of Mathematics and Statistics, Wuhan University,
\\Wuhan 430072, Hubei, People's Republic of China\\
$^{2}$ School of Mathematics, Sichuan University, Chengdu 610064, People's Republic of China \\
}
\date{}
\begin{document}
\newtheorem{definition}{Definition}[section]
\newtheorem{theorem}{Theorem}[section]
\newtheorem{lemma}{Lemma}[section]
\newtheorem{corollary}{Corollary}[section]
\newtheorem{example}{Example}[section]
\newtheorem{property}{Property}[section]
\newtheorem{proposition}{Proposition}[section]
\newtheorem{remark}{Remark}[section]

\newcommand{\qed}{\nolinebreak\hfill\rule{2mm}{2mm}
\par\medbreak}
\newcommand{\Proof}{\par\medbreak\it Proof: \rm}
\newcommand{\rem}{\par\medbreak\it Remark: \rm}
\newcommand{\defi}{\par\medbreak\it Definition : \rm}
\renewcommand{\thefootnote}{\arabic{footnote}}

\maketitle

\begin{abstract}
{\it In this paper, we establish first the resonance identity for non-contractible homologically visible prime closed geodesics on
Finsler $n$-dimensional real projective space $(\mathbb{R}P^n,F)$ when there exist only finitely many distinct non-contractible closed
geodesics on $(\mathbb{R}P^n,F)$, where the integer $n\geq2$. Then as an application of this resonance identity, we prove the existence
of at least two distinct non-contractible closed geodesics on $\RP$ with a bumpy and irreversible Finsler metric. Together with two previous
results on bumpy and reversible Finsler metrics in \cite{DLX2015} and \cite{Tai2016}, it yields that every $\RP$ with a bumpy Finsler metric
possesses at least two distinct non-contractible closed geodesics.}
\end{abstract}



\renewcommand{\theequation}{\thesection.\arabic{equation}}
\makeatletter\@addtoreset{equation}{section}\makeatother

\setcounter{equation}{0}
\section{Introduction}

In this paper, we are concerned with the multiplicity of closed geodesics on
$n$-dimensional real projective space $\mathbb{R}P^n$ with a Finsler metric $F$,
which is the typically non-simply connected manifold with the fundamental group $\Z_{2}$.
One of the main ingredients is a new resonance identity of non-contractible homologically visible prime closed
geodesics on $(\mathbb{R}P^n,F)$ when there exist only finitely many distinct non-contractible closed
geodesics on $(\mathbb{R}P^n,F)$. The second one is the precise iteration formulae of Morse indices for non-orientable
closed geodesics which can be seen as a complement of the index iteration theory for the orientable case. The third one
is the application of Kronecker's approximation theorem in Number theory to the multiplicity of non-contractible closed
geodesics on $(\mathbb{R}P^n,F)$.

A closed curve on a Finsler manifold is a closed geodesic if it is
locally the shortest path connecting any two nearby points on this
curve. As usual, on any Finsler manifold
$(M, F)$, a closed geodesic $c:S^1=\R/\Z\to M$ is {\it prime}
if it is not a multiple covering (i.e., iteration) of any other
closed geodesics. Here the $m$-th iteration $c^m$ of $c$ is defined
by $c^m(t)=c(mt)$. The inverse curve $c^{-1}$ of $c$ is defined by
$c^{-1}(t)=c(1-t)$ for $t\in \R$.  Note that unlike Riemannian manifold,
the inverse curve $c^{-1}$ of a closed geodesic $c$
on a irreversible Finsler manifold need not be a geodesic.
We call two prime closed geodesics
$c$ and $d$ {\it distinct} if there is no $\th\in (0,1)$ such that
$c(t)=d(t+\th)$ for all $t\in\R$.
For a closed geodesic $c$ on $(M,\,F)$, denote by $P_c$
the linearized Poincar\'{e} map of $c$. Recall that a Finsler metric $F$ is {\it bumpy} if all the closed geodesics
on $(M, \,F)$ are non-degenerate, i.e., $1\notin \sigma(P_c)$ for any closed
geodesic $c$.

Let $\Lm M$ be the free loop space on $M$ defined by
\begin{equation}
\label{LambdaM}
  \Lambda M=\left\{\gamma: S^{1}\to M\mid \gamma\ {\rm is\ absolutely\ continuous\ and}\
                        \int_{0}^{1}F(\gamma,\dot{\gamma})^{2}dt<+\infty\right\},
\end{equation}
endowed with a natural structure of Riemannian Hilbert manifold on which the group $S^1=\R/\Z$ acts continuously by
isometries (cf. Shen \cite{Shen2001}).

It is well known (cf. Chapter 1 of Klingenberg \cite{Kli1978}) that $c$ is a closed geodesic or a constant curve
on $(M,F)$ if and only if $c$ is a critical point of the energy functional
\begin{equation}
\label{energy}
E(\gamma)=\frac{1}{2}\int_{0}^{1}F(\gamma,\dot{\gamma})^{2}dt.
\end{equation}
Based on it, many important results on this subject have been obtained (cf. \cite{Ano}, \cite{Ban}, \cite{Fra}, \cite{Hin1984}-\cite{Hin1993}, \cite{Rad1989}-\cite{Rad1992}).
In particular, in 1969 Gromoll and Meyer \cite{GM1969JDG} used Morse theory and Bott's index iteration formulae \cite{Bott1956}
to establish the existence of infinitely many distinct closed geodesics on $M$, when the Betti number sequence
$\{{\beta}_k(\Lm M;\mathbb{Q})\}_{k\in\Z}$ is unbounded. Then Vigu$\acute{e}$-Poirrier and Sullivan \cite{VS1976}
further proved in 1976 that for a compact simply connected manifold $M$, the Gromoll-Meyer condition holds if and
only if $H^{*}(M;\mathbb{Q})$ is generated by more than one element. Here the Gromoll-Meyer theorem is valid actually
for any field $\mathbb{F}$. Note that it can not be applied to the compact rank one symmetric spaces
\be   S^{n},\ \mathbb{R}P^{n},\ \mathbb{C}P^{n},\ \mathbb{H}P^{n}\ \text{and}\ {\rm CaP}^{2},  \label{mflds}\ee
since $\{{\beta}_k(\Lm M,\mathbb{F})\}_{k\in\Z}$ with $M$ in (\ref{mflds}) is bounded with respect to any field
$\mathbb{F}$ (cf. Ziller \cite{Ziller1977}). In fact, each of them endowed with Katok metrics possesses only finitely many distinct prime closed geodesics
 (cf. Katok \cite{Katok1973}, also Ziller \cite{Ziller1982}).

In 2005,  Bangert and Long \cite{BL2010} (published in 2010) showed the existence of at least two distinct
closed geodesics on every Finsler $S^2$. Subsequently, such a multiplicity result for
$S^{n}$ with a bumpy Finsler metric, i.e., on which all closed geodesics are non-degenerate, was proved by Duan and Long
\cite{DL2007} and Rademacher \cite{Rad2010} independently.  In recent years, more interesting results on this problem have been obtained,
such as \cite{DL2010}-\cite{DLW2}, \cite{HiR}, \cite{LD2009}, \cite{Rad2007}, \cite{Wang2008}-\cite{Wang2012}.  We refer the readers to the survey
papers of Long \cite{lo2006}, Taimanov \cite{Tai2010} and Oancea \cite{Oancea2014} for more studies on this subject.

Besides many works on closed geodesics in the literature which study closed geodesics on simply connected
manifolds, we are aware of  not many papers on the multiplicity of closed geodesics on non-simply connected ones
published before 2015, at least when they are endowed with Finsler metrics. For example,
Ballman, Thorbergsson and Ziller \cite{BTZ1981} of 1981 and Bangert and Hingston \cite{BH1984} of 1984 dealt with the
non-simply connected manifolds with a finite/infinite cyclic fundamental group respectively by the min-max principle.

In order to apply Morse theory to the multiplicity of closed geodesics on $\RP$, motivated by the studies on the simply
connected manifolds, in particular, the resonance identity proved by Rademacher \cite{Rad1989}, Xiao and Long \cite{XL2015}
in 2015 investigated the topological structure of the non-contractible
loop space and established the resonance identity for the non-contractible closed geodesics on $\R P^{2n+1}$ by using
$\Z_2$ coefficient homology. As an application, Duan, Long and Xiao \cite{DLX2015} proved the existence of at least two
distinct non-contractible closed geodesics on $\R P^{3}$ endowed with a bumpy and irreversible Finsler metric. In a very recent paper \cite{Tai2016}, Taimanov
studied the rational equivariant cohomology of the spaces of non-contractible loops in compact space forms and
proved the existence of at least two distinct non-contractible closed geodesics on $\mathbb{R}P^2$ endowed with a bumpy and irreversible Finsler metric.
Then Liu \cite{Liuhui2016} combined Fadell-Rabinowitz index theory with Taimanov's topological results to get multiplicity
results of non-contractible closed geodesics on positively curved Finsler $\RP$.

Motivated by \cite{Tai2016} and \cite{XL2015}, in section 2 of this paper we obtain the resonance identity for the non-contractible closed geodesics on $\R P^{n}$
by using rational coefficient homology for any $n\geq2$ regardless of whether $n$ is odd or not.
\begin{theorem}\label{Thm1.1} Suppose the Finsler manifold $M=(\mathbb{R}P^{n},F)$ possesses only finitely
many distinct non-contractible prime closed geodesics, among which we denote the distinct non-contractible homologically visible prime
closed geodesics by $c_1, \ldots, c_r$  for some integer $r>0$, where $n\geq2$. Then we have
\bea  \sum_{j=1}^{r}\frac{\hat{\chi}(c_j)}{\hat{i}(c_j)} = \bar{B}(\Lm_gM ;\Q) =\left\{\begin{array}{ll}
    \frac{n+1}{2(n-1)},&\ if\ n\in 2\N-1,\\
    \frac{n}{2(n-1)},&\ if\ n\in 2\N.\\
    \end{array}\right.     \lb{reident1}\eea
where $\Lm_{g} M$ is the non-contractible loop space of $M$ and the mean Euler number $\hat{\chi}(c_j)$ of $c_j$ is defined by
$$  \hat{\chi}(c_j) = \frac{1}{n_j}\sum_{m=1}^{n_j/2}\sum_{l=0}^{2n-2}(-1)^{l+i(c_{j})}k_{l}(c_{j}^{2m-1})\in\Q, $$
and $n_j=n_{c_j}$ is the analytical period of $c_j$, $k_{l}(c_{j}^{2m-1})$ is the local homological type number of $c_{j}^{2m-1}$,
$i(c_{j})$ and $\hat{i}(c_j)$ are the Morse index and mean index of $c_j$ respectively.

In particular, if the Finsler metric $F$ on $\mathbb{R}P^{n}$ is bumpy, then (\ref{reident1}) has the following simple form
\bea \sum_{j=1}^{r}\frac{(-1)^{i(c_j)}}{\hat{i}(c_j)} =\left\{\begin{array}{ll}
    \frac{n+1}{n-1},&\ if\ n\in 2\N-1,\\
    \frac{n}{n-1},&\ if\ n\in 2\N.\\
    \end{array}\right.     \lb{breident1}\eea
\end{theorem}

Based on Theorem \ref{Thm1.1}, the precise iteration formulae of Morse indices for closed geodesics and Morse theory,  especially the $S^1$-equivariant Poincar$\acute{e}$ series of $\Lm_{g}M$ derived by Taimanov (cf. Lemma \ref{Lm2.3}), and using some techniques in Number theory, we can prove the following multiplicity result of non-contractible closed geodesics on $(\mathbb{R}P^{n},F)$.

\begin{theorem}
\label{mainresult}
Every $\RP$ endowed with a bumpy and irreversible Finsler metric $F$ has at least two distinct non-contractible closed geodesics, where $n\geq2$.
\end{theorem}
\begin{remark}
 For any compact simply-connected bumpy Finsler manifold,  Duan, Long and Wang in \cite{DLW1} proved the same conclusion as Theorem \ref{mainresult}. However, their method is not applicable to our problem. Indeed, one of the crucial facts in their proof is that if there is only one prime closed geodesic on such a manifold, its Morse index must be greater than or equal to some positive integer. But there is always a minimal point of the energy functional on $\Lm_{g}(\R P^{n})$ with Morse index 0.
 \end{remark}

If $F$ is a bumpy and reversible Finsler metric, the same conclusion as Theorem \ref{mainresult} has been proved in Theorem 1.2 of \cite{DLX2015}
and the remark behind Theorem 5 of \cite{Tai2016}.  As a combined outcome, we immediately get the desired result as follows.
  \begin{corollary}
  Every $\RP$ endowed with a bumpy Finsler metric has at least two distinct non-contractible closed geodesics, where $n\geq2$.
   \end{corollary}

This paper is organized as follows. In section 2, we use Morse theory to establish the resonance identity of Theorem 1.1. Then in section 3,
we investigate the precise iteration formulae of Morse indices for closed geodesics on $\R P^{n}$ and build a bridge between their Morse indices  and a division of an interval.  In section 4, a special system of irrational numbers associated to our problem is carefully studied and  a key result on it for our later proof of Theorem \ref{mainresult} is obtained. Finally in section 5, we draw support from the well known Kronecker's
approximation theorem in Number theory and give the proof of Theorem \ref{mainresult}.

We close this introduction with some illustrations of notations in this paper. As usual, let $\N$,  $\Z$, $\Q$ and $\Q^{c}$ denote the sets of natural integers,
integers, rational numbers and  irrational numbers respectively. We also use  notations $E(a)=\min\{k\in\Z\,|\,k\ge a\}$,
$[a]=\max\{k\in\Z\,|\,k\le a\}$, $\varphi(a)=E(a)-[a]$ and $\{a\}=a-[a]$ for any $a\in\R$. Throughout this paper, we use $\Q$ coefficients for all homological and cohomological modules.

\setcounter{equation}{0}
\section{Morse theory and resonance identity of non-contractible closed geodesics on $(\mathbb{R}P^n,F)$}
Let $M=(M,F)$ be a compact Finsler manifold, the space
$\Lambda=\Lambda M$ of $H^1$-maps $\gamma:S^1\rightarrow M$ has a
natural structure of Riemannian Hilbert manifolds on which the
group $S^1=\R/\Z$ acts continuously by isometries. This action is defined by
$(s\cdot\gamma)(t)=\gamma(t+s)$ for all $\gamma\in\Lm$ and $s,
t\in S^1$. For any $\gamma\in\Lambda$, the energy functional is
defined by
\be E(\gamma)=\frac{1}{2}\int_{S^1}F(\gamma(t),\dot{\gamma}(t))^2dt.
\lb{2.1}\ee
It is $C^{1,1}$ and invariant under the $S^1$-action. The
critical points of $E$ of positive energies are precisely the closed geodesics
$\gamma:S^1\to M$. The index form of the functional $E$ is well
defined along any closed geodesic $c$ on $M$, which we denote by
$E''(c)$. As usual, we denote by $i(c)$ and
$\nu(c)-1$ the Morse index and nullity of $E$ at $c$. In the
following, we denote by
\be \Lm^\kappa=\{d\in \Lm\;|\;E(d)\le\kappa\},\quad \Lm^{\kappa-}=\{d\in \Lm\;|\; E(d)<\kappa\},
  \quad \forall \kappa\ge 0. \lb{2.2}\ee
For a closed geodesic $c$ we set $ \Lm(c)=\{\ga\in\Lm\mid E(\ga)<E(c)\}$.

For $m\in\N$ we denote the $m$-fold iteration map
$\phi_m:\Lambda\rightarrow\Lambda$ by $\phi_m(\ga)(t)=\ga(mt)$, for all
$\,\ga\in\Lm, t\in S^1$, as well as $\ga^m=\phi_m(\gamma)$. If $\gamma\in\Lambda$
is not constant then the multiplicity $m(\gamma)$ of $\gamma$ is the order of the
isotropy group $\{s\in S^1\mid s\cdot\gamma=\gamma\}$. For a closed geodesic $c$,
the mean index $\hat{i}(c)$ is defined as usual by
$\hat{i}(c)=\lim_{m\to\infty}i(c^m)/m$. Using singular homology with rational
coefficients we consider the following critical $\Q$-module of a closed geodesic
$c\in\Lambda$:
\be \overline{C}_*(E,c)
   = H_*\left((\Lm(c)\cup S^1\cdot c)/S^1,\Lm(c)/S^1; \Q\right). \lb{2.3}\ee
In the following we let $M=\mathbb{R}P^{n}$, where $n\geq2$, it is well known that $\pi_1(\mathbb{R}P^{n})=\Z_2=\{e, g\}$ with $e$ being
the identity and $g$ being the generator of $\Z_2$ satisfying $g^2=e$. Then the free loop space $\Lambda M$ possesses
a natural decomposition\bea \Lambda M=\Lambda_e M\bigsqcup\Lambda_g M,\nn\eea
where $\Lambda_e M$ and $\Lambda_g M$ are the two connected components of $\Lambda M$ whose elements are homotopic
to $e$ and $g$ respectively. We set $\Lambda_{g}(c) = \{\gamma\in \Lambda_{g}M\mid E(\gamma)<E(c)\}$.
Note that for a non-contractible prime closed geodesic $c$, $c^m\in\Lambda_g M$ if and only if $m$ is odd.

We call a non-contractible prime closed geodesic satisfying the isolation condition, if
the following holds:

{\bf (Iso)  For all $m\in\N$ the orbit $S^1\cdot c^{2m-1}$ is an
isolated critical orbit of $E$. }

Note that if the number of non-contractible prime closed geodesics on $M=\mathbb{R}P^{n}$
is finite, then all the non-contractible prime closed geodesics satisfy (Iso).

If a non-contractible closed geodesic $c$ has multiplicity $2m-1$, then the subgroup $\Z_{2m-1}=\{\frac{l}{2m-1}\mid 0\leq l<2m-1\}$
of $S^1$ acts on $\overline{C}_*(E,c)$. As studied in p.59 of \cite{Rad1992},
for all $m\in\N$, let
$H_{\ast}(X,A)^{\pm\Z_{2m-1}}
   = \{[\xi]\in H_{\ast}(X,A)\,|\,T_{\ast}[\xi]=\pm [\xi]\}$,
where $T$ is a generator of the $\Z_{2m-1}$-action.
On $S^1$-critical modules of $c^{2m-1}$, the following lemma holds:
\begin{lemma}
\label{Rad1992} {\rm (cf. Satz 6.11 of \cite{Rad1992} and \cite{BL2010})}  Suppose $c$ is
a non-contractible prime closed geodesic on a Finsler manifold $M=\mathbb{R}P^{n}$ satisfying (Iso). Then
there exist $U_{c^{2m-1}}$ and $N_{c^{2m-1}}$, the so-called local negative
disk and the local characteristic manifold at $c^{2m-1}$ respectively,
such that $\nu(c^{2m-1})=\dim N_{c^{2m-1}}$ and
\bea &&\overline{C}_q( E,c^{2m-1})
\equiv H_q\left((\Lm_g(c^{2m-1})\cup S^1\cdot c^{2m-1})/S^1, \Lm_g(c^{2m-1})/S^1\right)\nn\\
&=& \left(H_{i(c^{2m-1})}(U_{c^{2m-1}}^-\cup\{c^{2m-1}\},U_{c^{2m-1}}^-)
    \otimes H_{q-i(c^{2m-1})}(N_{c^{2m-1}}^-\cup\{c^{2m-1}\},N_{c^{2m-1}}^-)\right)^{+\Z_{2m-1}}, \nn
\eea
where $U_{c^{2m-1}}^-=U_{c^{2m-1}}\cap\Lm_g(c^{2m-1})$ and $N_{c^{2m-1}}^-=N_{c^{2m-1}}\cap\Lm_g(c^{2m-1})$.

(i) When $\nu(c^{2m-1})=0$, there holds
\bea \overline{C}_q( E,c^{2m-1}) = \left\{\begin{array}{ll}

     \Q, &\quad {\it if}\;\;
                   q=i(c^{2m-1}),\;  \cr
     0, &\quad {\it otherwise},\\ \end{array}\right. \nn \eea

(ii) When $\nu(c^{2m-1})>0$, there holds
$$ \overline{C}_q( E,c^{2m-1})=H_{q-i(c^{2m-1})}(N_{c^{2m-1}}^-\cup\{c^{2m-1}\},N_{c^{2m-1}}^-)^{+\Z_{2m-1}}, $$
where we have used the fact $i(c^{2m-1})-i(c)\in 2\Z$.
\end{lemma}

As usual, for $m\in\N$ and $l\in\Z$ we define the local homological type numbers of $c^{2m-1}$ by
\be k_{l}(c^{2m-1})
= \dim H_{l}(N_{c^{2m-1}}^{-}\cup\{c^{2m-1}\},N_{c^{2m-1}}^{-})^{+\Z_{2m-1}}.  \lb{CGht1}\ee

Based on works of Rademacher in \cite{Rad1989}, Long and Duan in \cite{LD2009} and \cite{DL2010},
we define the {\it analytical period} $n_c$ of the closed geodesic $c$ by
\be n_c = \min\{j\in 2\N\,|\,\nu(c^j)=\max_{m\ge 1}\nu(c^m),\;\;
                  \forall\,m\in 2\N-1\}. \lb{CGap1}\ee
Note that here in order to simplify the study for non-contractible closed geodesics in $\mathbb{R}P^{n}$,
we have slightly modified the definition in \cite{LD2009} and \cite{DL2010} by requiring the analytical
period to be even. Then by the same proofs in \cite{LD2009} and \cite{DL2010}, we have
\be  k_{l}(c^{2m-1+hn_c}) = k_{l}(c^{2m-1}), \qquad \forall\;m,\;h\in \N,\;l\in\Z.  \lb{CGap2}\ee
For more detailed properties of the analytical period $n_c$ of a closed geodesic $c$, we refer readers to
the two Section 3s in \cite{LD2009} and \cite{DL2010}.

As in \cite{BK1983}, we have
\begin{definition}\label{def-hv} Let $(M,F)$ be a compact Finsler manifold. A closed geodesic $c$ on $M$
is homologically visible, if there exists an integer $k\in\Z$ such that $\bar{C}_k(E,c) \not= 0$. We denote
by $\CG_{\hv}(M,F)$ the set of all distinct homologically visible prime closed geodesics on $(M,F)$.
\end{definition}
\begin{lemma}\label{Lm2.2} Suppose the Finsler manifold $M=(\mathbb{R}P^{n},F)$ possesses only finitely
many distinct non-contractible prime closed geodesics, among which we denote the distinct non-contractible homologically visible prime closed
geodesics by $c_1, \ldots, c_r$  for some integer $r>0$. Then we have
\be  \hat{i}(c_i) > 0, \qquad \forall\;1\le i\le r.  \lb{fcg.1}\ee
\end{lemma}

\Proof First, we claim that Theorem 3 in \cite{BK1983} for $M=\RP$ can be stated as:

`` Let $c$ be a closed geodesic in $\Lambda_{g}M$ such that $i(c^{m})=0$ for all $m\in\N$. Suppose $c$ is neither
homologically invisible nor an absolute minimum of $E$ in $\Lambda_{g}M$. Then there exist infinitely many closed geodesics in $\Lambda_{g}M$."

Indeed, one can focus the proofs of Theorem 3 in \cite{BK1983} on $\Lambda_{g}M$ with some obvious modifications.
Assume by contradiction. Similarly as in \cite{BK1983}, we can choose a different $c\in\Lm_{g}M$, if necessary, and find $p\in\N$ such that
$H_{p}(\Lambda_{g}(c)\cup S\cdot c,\Lambda_{g}(c))\neq0$ and $H_{q}(\Lambda_{g}(c)\cup S\cdot c,\Lambda_{g}(c))=0$ for every $q>p$ and
every closed geodesic $d\in\Lm_{g}M$ with $i(d^{m})\equiv0.$

Consider the following commutative diagram \be
\begin{tabular}{ccc}
$H_{p}(\Lambda_{g}(c)\cup S\cdot c,\Lambda_{g}(c))$   & $\mapright{{\psi}^{m}_{*}}$ & $H_{p}(\Lambda_{g}(c^{m})\cup S\cdot c^{m},\Lambda_{g}(c^{m}))$ \\
$\mapdown{i_*}$&                       & $\mapdown{i_*}$ \\
$H_{p}(\Lambda_{g}M,\Lambda_{g}(c))$                     & $\mapright{{\psi}^{m}_{*}}$ & $H_{p}(\Lambda_{g}M,\Lambda_{g}(c^{m}))$,\\
\end{tabular}  \lb{926diagram}\ee
where $m$ is odd and $\psi^{m}:\Lm_{g} M\to\Lm_{g} M$ is the $m$-fold iteration map. By similar arguments as those in \cite{BK1983},
there is $A>0$ such that the map $i_{*}\circ\psi^{m}_{*}$ is one-to-one, if $E(c^{m})>A$ and none of the $k_{i}\in K_0$ divides $m$
where  $$K_{0}=\{k_{0},k_{1},k_{2},\dots,k_{s}\},$$ with $k_{0}=2$ and $k_{1},k_{2},\dots,k_{s}$ therein. Here note that the
required $m$ is odd and so $c^{m}\in\Lambda_{g}(M)$ for $c\in\Lambda_{g}M$.

On the other hand, we define $$K=\{m\geq2\mid E(c^{m})\leq A\}\cup K_{0}.$$
Then by Corollary 1 of \cite{BK1983}, there exists $\bar{m}\in\N\backslash\{1\}$ such that no $k\in K$ divides $\bar{m}$ and
$\psi^{\bar{m}}_{*}\circ i_{*}$ vanishes.  In particular, $E(c^{\bar{m}})>A$ and none of the $k_{i}\in K_0$ divides $\bar{m}$.
Due to $\psi^{\bar{m}}_{*}\circ i_{*}=i_{*}\circ \psi^{\bar{m}}_{*}$ in (\ref{926diagram}),  this yields a contradiction.
Hence there exist infinitely many closed geodesics in $\Lambda_{g}M$.

Accordingly, Corollary 2 in \cite{BK1983}  for $M=\RP$ can be stated as:

``  Suppose there exists a closed geodesic $c\in\Lambda_{g}M$ such that $c^{m}$ is a local minimum of $E$ in $\Lambda_{g}M$ for
infinitely many odd $m\in\N$. Then there exist infinitely many closed geodesics in $\Lambda_{g}M$."

Based on the above two variants of Theorem 3 and Corollary 2 in \cite{BK1983}, we can prove our Lemma \ref{Lm2.2} as follows.

It is well known that every closed geodesic $c$ on $M$ must have mean index $\hat{i}(c)\ge 0$.

Assume by contradiction that there is a non-contractible homologically visible prime closed geodesic $c$ on $M$ satisfying
$\hat{i}(c)=0$. Then $i(c^m)=0$ for all $m\in\N$ by Bott iteration formula and $c$ must be an absolute
minimum of $E$ in $\Lambda_{g}M$, since otherwise there would exist infinitely many distinct non-contractible
closed geodesics on $M$ by the above variant of Theorem 3 on p.385 of \cite{BK1983}.

On the other hand, by Lemma 7.1 of \cite{Rad1992}, there exists a $k(c)\in2\N$ such that
$\nu(c^{m+k(c)})=\nu(c^m)$ for all $m\in\N$. Specially we obtain $\nu(c^{mk(c)+1}) = \nu(c)$ for all
$m\in\N$ and then elements of $\ker(E''(c^{mk(c)+1}))$ are precisely $mk(c)+1$st iterates of
elements of $\ker(E''(c))$. Thus by the Gromoll-Meyer theorem in \cite{GM1969Top}, the behavior of the
restriction of $E$ to $\ker(E''(c^{mk(c)+1}))$ is the same as that of the restriction of $E$ to
$\ker(E''(c))$. Then together with the fact $i(c^m)=0$ for all $m\in\N$, we obtain that $c^{mk(c)+1}$
is a local minimum of $E$ in $\Lambda_{g}M$ for every $m\in\N$. Because $M$ is compact and
possessing finite fundamental group, there must exist infinitely many distinct non-contractible closed geodesics on $M$
by the above variant of Corollary 2 on p.386 of \cite{BK1983}. Then it yields a contradiction and proves (\ref{fcg.1}).
$\hfill\Box$

In \cite{Tai2016}, Taimanov calculated the rational equivariant cohomology of the spaces of non-contractible loops of $\mathbb{R}P^n$
which is crucial for us to prove Theorem \ref{Thm1.1} and  can be stated as follows.
\begin{lemma}\label{Lm2.3} {\rm (cf. Theorem 3 of \cite{Tai2016} or Lemma 2.2 of \cite{Liuhui2016})} For $M=\mathbb{R}P^n$, we have

(i) When $n=2k+1$ is odd, the $S^1$-cohomology ring of $\Lambda_g M$ has the form
$$H^{S^1, *}(\Lambda_g M; \Q)=\Q[w, z]/ \{w^{k+1} = 0\}, \quad deg(w)=2, \quad deg(z)=2k$$
Then the $S^1$-equivariant Poincar$\acute{e}$ series
of $\Lambda_g M$ is given by
\bea P^{S^1}(\Lambda_g M; \Q)(t)&=&\frac{1-t^{2k+2}}{(1-t^2)(1-t^{2k})}\nn\\
&=&\frac{1}{1-t^2}+\frac{t^{2k}}{1-t^{2k}}\nn\\&=&(1+t^2+t^4+\cdots+t^{2k}+\cdots)+(t^{2k}+t^{4k}+t^{6k}+\cdots),\nn\eea
which yields Betti numbers
\bea \bar{\beta}_q=\rank H_q^{S^1}(\Lambda_g M;\Q)=\left\{\begin{array}{ll}
    2,&\ if\ q\in \{j(n-1)\mid j\in\mathbb{N}\},\\
    1,&\ if\ q\in (2\mathbb{N}\cup\{0\})\backslash\{j(n-1)\mid j\in\mathbb{N}\},\\
    0,&\ otherwise.\\
    \end{array}\right. \lb{b.1}\eea
and the average $S^1$-equivariant Betti number of $\Lambda_g M$ satisfies
\be  \bar{B}(\Lambda_g M;\Q)\equiv\lim_{q\to+\infty}\frac{1}{q}\sum_{k=0}^{q}(-1)^k\bar{\beta}_k= \frac{n+1}{2(n-1)}.   \lb{aB.1}\ee

(ii) When $n=2k$ is even, the $S^1$-cohomology ring of $\Lambda_g M$ has the form
$$H^{S^1, *}(\Lambda_g M; \Q)=\Q[w, z]/ \{w^{2k} = 0\}, \quad deg(w)=2, \quad deg(z)=4k-2$$
Then the $S^1$-equivariant Poincar$\acute{e}$ series
of $\Lambda_g M$ is given by \bea P^{S^1}(\Lambda_g M; \Q)(t)&=&\frac{1-t^{4k}}{(1-t^2)(1-t^{4k-2})}\nn\\
&=&\frac{1}{1-t^2}+\frac{t^{4k-2}}{1-t^{4k-2}}\nn\\&=&(1+t^2+t^4+\cdots+t^{2k}+\cdots)+(t^{4k-2}+t^{2(4k-2)}+t^{3(4k-2)}+\cdots),\nn\eea
which yields Betti numbers \bea \bar{\beta}_q= \rank H_q^{S^1}(\Lambda_g M;\Q)=\left\{\begin{array}{ll}
2,&\ if\ q\in \{2j(n-1)\mid j\in\mathbb{N}\},\\
1,&\ if\ q\in (2\mathbb{N}\cup\{0\})\backslash\{2j(n-1)\mid j\in\mathbb{N}\},\\
0,&\ otherwise.\\
\end{array}
\right.\lb{b.2}\eea
and the average $S^1$-equivariant Betti number of $\Lambda_g M$ satisfies
\be  \bar{B}(\Lambda_g M;\Q)\equiv\lim_{q\to+\infty}\frac{1}{q}\sum_{k=0}^{q}(-1)^k\bar{\beta}_k= \frac{n}{2(n-1)}.   \lb{aB.2}\ee
\end{lemma}
\begin{remark}
For the case of $\R P^{2n+1}$, the same conclusions as (\ref{b.1}) and (\ref{aB.1}) were obtained in \cite{XL2015} where the coefficient field $\Z_{2}$ was used and they are also effective to our problem  since the multiplicity of every curve on  $\Lambda_g M$ is odd.
\end{remark}

Now we give the proof of the resonance identity in Theorem \ref{Thm1.1}.

{\bf Proof of Theorem \ref{Thm1.1}.} Recall that we denote the non-contractible homologically visible prime closed
geodesics by $\CG_{\hv}(M)=\{c_1, \ldots, c_r\}$ for some integer $r>0$ when the number of distinct non-contractible prime
closed geodesics on $M=\mathbb{R}P^n$ is finite. Note also that by Lemma \ref{Lm2.2} we have $\hat{i}(c_j)>0$ for all
$1\le j\le r$.

Let \bea m_q\equiv M_q(\Lm_gM) =\sum_{1\le j\le r,\; m\ge 1}\dim{\ol{C}}_q(E, c^{2m-1}_j),\quad q\in\Z.\nn\eea
The Morse series of $\Lm_gM$ is defined by
\be  M(t) = \sum_{h=0}^{+\infty}m_ht^h.  \label{wh}\ee

{\bf Claim 1.} {\it $\{m_h\}$ is a bounded sequence.}

In fact, by (\ref{CGap2}), we have \bea m_h=\sum_{j=1}^{r}\sum_{m=1}^{n_j/2}\sum_{l=0}^{2n-2}k_{l}(c_{j}^{2m-1})
            \;{}^{\#}\left\{s\in\mathbb{N}\cup\{0\}\mid h-i(c_{j}^{2m-1+sn_j})=l\right\}, \label{wh2}\eea
and by Theorem 9.2.1, Theorems 10.1.2 of \cite{lo2000}, and Lemmas \ref{orient}-\ref{unorient} below,
we have $|i(c_{j}^{2m-1+sn_j})-(2m-1+sn_j)\hat{i}(c_{j})|\le 2n-2$, then
\bea
&&{}^{\#}\left\{s\in\mathbb{N}\cup\{0\}\mid h-i(c_{j}^{2m-1+sn_j})=l\right\} \nn\\
&&\qquad = \;{}^{\#}\left\{s\in\mathbb{N}\cup\{0\}\mid l+i(c_{j}^{2m-1+sn_j})=h,\;
           |i(c_{j}^{2m-1+sn_j})-(2m-1+sn_j)\hat{i}(c_{j})|\le 2n-2\right\}  \nn\\
&&\qquad \le \;{}^{\#}\left\{s\in\mathbb{N}\cup\{0\}\mid 2n-2\geq |h-l-(2m-1+sn_{j})\hat{i}(c_{j})|  \right\}  \nn\\
&&\qquad = \;{}^{\#}\left\{s\in\mathbb{N}\cup\{0\}\mid \frac{h-l-2n+2-(2m-1)\hat{i}(c_j)}{n_j\hat{i}(c_j)} \leq s
       \le \frac{h-l+2n-2-(2m-1)\hat{i}(c_j)}{n_j\hat{i}(c_j)}\right\} \nn\\
&&\qquad \le \frac{4n-4}{n_j\hat{i}(c_j)} + 1.  \label{wh3}\eea
Hence Claim 1 follows by (\ref{wh2}) and (\ref{wh3}).

We now use the method
in the proof of Theorem 5.4 of \cite{LW2007} to estimate
$$   M^{q}(-1) = \sum_{h=0}^{q}m_h(-1)^h.   $$

By (\ref{wh}) and (\ref{CGap2}) we obtain
\bea M^{q}(-1)
&=& \sum_{h=0}^{q}m_{h}(-1)^{h}  \nn\\
&=& \sum_{j=1}^{r}\sum_{m=1}^{n_j/2}\sum_{l=0}^{2n-2}\sum_{h=0}^{q}(-1)^{h}k_{l}(c_{j}^{2m-1})
            \;{}^{\#}\left\{s\in\mathbb{N}\cup\{0\}\mid h-i(c_{j}^{2m-1+sn_j})=l\right\}  \nn\\
&=& \sum_{j=1}^{r}\sum_{m=1}^{n_j/2}\sum_{l=0}^{2n-2}(-1)^{l+i(c_{j})}k_{l}(c_{j}^{2m-1})
            \;{}^{\#}\left\{s\in\mathbb{N}\cup\{0\}\mid l+i(c_{j}^{2m-1+sn_j})\le q\right\}.  \nn\eea
On the one hand, we have
\bea
&&{}^{\#}\left\{s\in\mathbb{N}\cup\{0\}\mid l+i(c_{j}^{2m-1+sn_j})\le q\right\} \nn\\
&&\qquad = \;{}^{\#}\left\{s\in\mathbb{N}\cup\{0\}\mid l+i(c_{j}^{2m-1+sn_j})\le q,\;
           |i(c_{j}^{2m-1+sn_j})-(2m-1+sn_j)\hat{i}(c_{j})|\le 2n-2\right\}  \nn\\
&&\qquad \le \;{}^{\#}\left\{s\in\mathbb{N}\cup\{0\}\mid 0\le (2m-1+sn_{j})\hat{i}(c_{j})\le q-l+2n-2  \right\}  \nn\\
&&\qquad = \;{}^{\#}\left\{s\in\mathbb{N}\cup\{0\}\mid 0 \leq s
       \le \frac{q-l+2n-2-(2m-1)\hat{i}(c_j)}{n_j\hat{i}(c_j)}\right\} \nn\\
&&\qquad \le \frac{q-l+2n-2}{n_j\hat{i}(c_j)} + 1. \nn\eea
On the other hand, we have
\bea
&&{}^{\#}\left\{s\in\mathbb{N}\cup\{0\}\mid l+i(c_{j}^{2m-1+sn_j})\le q \right\} \nn\\
&&\qquad = \;{}^{\#}\left\{s\in\mathbb{N}\cup\{0\}\mid l+i(c_{j}^{2m-1+sn_j})\le q,\;
           |i(c_{j}^{2m-1+sn_j})-(2m-1+sn_{j})\hat{i}(c_{j})|\le 2n-2\right\}  \nn\\
&&\qquad \ge \;{}^{\#}\left\{s\in\mathbb{N}\cup\{0\}\mid i(c_{j}^{2m-1+sn_j})
       \le (2m-1+sn_{j})\hat{i}(c_{j})+2n-2\le q-l \right\} \nn\\
&&\qquad \ge \;{}^{\#}\left\{s\in\mathbb{N}\cup\{0\}\mid 0 \le s
       \le \frac{q-l-2n+2-(2m-1)\hat{i}(c_{j})}{n_j\hat{i}(c_{j})} \right\}  \nn\\
&&\qquad \ge \frac{q-l-2n+2}{n_j\hat{i}(c_{j})} - 1.  \nn\eea
Thus we obtain
$$  \lim_{q\to+\infty}\frac{1}{q}M^{q}(-1)
  = \sum_{j=1}^{r}\sum_{m=1}^{n_j/2}\sum_{l=0}^{2n-2}(-1)^{l+i(c_{j})}k_{l}(c_{j}^{2m-1})\frac{1}{n_{j}\hat{i}(c_{j})}
              = \sum_{j=1}^{r}\frac{\hat{\chi}(c_j)}{\hat{i}(c_j)}.  $$
Since $m_{h}$ is bounded, we then obtain
$$  \lim_{q\to+\infty}\frac{1}{q}M^{q}(-1) = \lim_{q\to+\infty}\frac{1}{q}P^{S^{1},q}(\Lambda_{g}M;\Q)(-1)
          = \lim_{q\to+\infty}\frac{1}{q}\sum_{k=0}^{q}(-1)^k\bar{\beta}_k = \bar{B}(\Lm_gM;\Q),  $$
where $P^{S^{1},q}(\Lambda_{g}M;\Q)(t)$ is the truncated polynomial of
$P^{S^{1}}(\Lambda_{g}M;\Q)(t)$ with terms of degree less than or equal to $q$.
Thus by (\ref{aB.1}) and (\ref{aB.2}) we get
$$  \sum_{j=1}^{r}\frac{\hat{\chi}(c_j)}{\hat{i}(c_j)} = \left\{\begin{array}{ll}
    \frac{n+1}{2(n-1)},&\ if\ n\in 2\N-1,\\
    \frac{n}{2(n-1)},&\ if\ n\in 2\N.\\
    \end{array}\right. $$
which proves (\ref{reident1}) of Theorem \ref{Thm1.1}. For the special case when each $c_{j}^{2m-1}$ is non-degenerate with $1\leq j\leq r$ and $m\in\mathbb{N}$,
we have $n_{j}=2$ and $k_l(c_j)=1$ when $l=0$, and $k_l(c_j)=0$ for all other
$l\in\Z$. Then (\ref{reident1}) has the following simple form
\be  \sum_{j=1}^{r}(-1)^{i(c_{j})}\frac{1}{\hat{i}(c_{j})}=\left\{\begin{array}{ll}
    \frac{n+1}{n-1},&\ if\ n\in 2\N-1,\\
    \frac{n}{n-1},&\ if\ n\in 2\N.\\
    \end{array}\right.  \nn\ee
which proves (\ref{breident1}) of Theorem \ref{Thm1.1}. $\hfill\Box$

\setcounter{equation}{0}
\section{Index  iteration theory for closed geodesics}
\subsection{Index iteration formulae  for closed geodesics}
In \cite{lo1999} of 1999, Y. Long established the basic normal form
decomposition of symplectic matrices. Based on it, he
further established the precise iteration formulae of Maslov $\omega$-indices for
symplectic paths in \cite{lo2000}, which can be related to  Morse indices  of either orientable or non-orientable closed geodesics in a slightly different way (cf. \cite{Liu2005} and Chap. 12 of \cite{lo2002}). Roughly speaking, the orientable (resp. non-orientable) case  corresponds to $i_{1}$ (resp. $i_{-1}$) index,  where $i_1$ and $i_{-1}$ denote the cases of $\omega$-index with $\omega=1$ and $\omega=-1$ respectively (cf. Chap. 5 of \cite{lo2002}).  Although we are concerned with $\R P^{n}$ in this paper, we will state such a relation precisely in a general form due to its independent interest. Throughout this section we denote the Morse index of a closed geodesic $c$ by  ind($c$) in stead of $i(c)$ to avoid confusion of notations and write $i_1(\gamma)$  as $i(\gamma)$ for short.

For the reader's convenience, we  briefly review some basic materials  in Long's book \cite{lo2002}.

Let $P$ be a symplectic matrix in Sp$(2N-2)$ and $\Omega^{0}(P)$ be the path connected component of its homotopy set $\Omega(P)$ which contains $P$. Then there is a path $f\in C([0,1],\Omega^{0}(P))$
such that $f(0)=P$ and
 \bea f(1)
&=& N_1(1,1)^{\dm p_-}\,\dm\,I_{2p_0}\,\dm\,N_1(1,-1)^{\dm p_+}\nn\\
& &  \dm\,N_1(-1,1)^{\dm q_-}\,\dm\,(-I_{2q_0})\,\dm\,N_1(-1,-1)^{\dm q_+} \nn\\
&&\dm\,R(\th_1)\,\dm\,\cdots\,\dm\,R(\th_{r'})\,\dm\,R(\th_{r'+1})\,\dm\,\cdots\,\dm\,R(\th_r)\lb{nmf}\\
&&\dm\,N_2(e^{i\aa_{1}},A_{1})\,\dm\,\cdots\,\dm\,N_2(e^{i\aa_{r_{\ast}}},A_{r_{\ast}})\nn\\
& &  \dm\,N_2(e^{i\bb_{1}},B_{1})\,\dm\,\cdots\,\dm\,N_2(e^{i\bb_{r_{0}}},B_{r_{0}})\nn\\
& &\dm\,H(\pm 2)^{\dm h},\nn\eea
where $N_{1}(\lambda,\chi)=\left(
  \begin{array}{ll}
    \lambda\quad \chi\\
   0\quad \lambda\\
  \end{array}
  \right)$ with $\lambda=\pm 1$ and $\chi=0,\ \pm 1$; $H(b)=\left(
  \begin{array}{ll}
    b\quad 0\\
   0\quad b^{-1}\\
  \end{array}
  \right)$ with $b=\pm 2$;
  $R(\theta)= \left(
  \begin{array}{ll}
    \cos\theta\ -\sin\theta\\
   \sin\theta\quad\ \cos\theta\\
  \end{array}
  \right)$ with $\theta\in (0,2\pi)\setminus\{\pi\}$ and we suppose that $\pi<\theta_{j}<2\pi$ iff $1\leq j\leq
r'$;  $$N_{2}(e^{i\aa_{j}},A_{j})=\left(
  \begin{array}{ll}
    R(\aa_{j})\  A_{j}\\
  \ 0\quad\ R(\aa_{j})\\
  \end{array}
  \right)\ \text{and}\ N_{2}(e^{i\bb_{j}},B_{j})=\left(\begin{array}{ll}
    R(\bb_{j})\  B_{j}\\
  \ 0\quad\ R(\bb_{j})\\
  \end{array}
  \right)$$ with $\alpha_{j},\beta_{j}\in (0,2\pi)\setminus\{\pi\}$ are non-trivial and trivial basic normal forms respectively.

  Let $\gamma_{0}$ and $\gamma_{1}$ be two symplectic paths in Sp$(2N-2)$ connecting the identity matrix $I$ to $P$ and $f(1)$
  satisfying $\gamma_{0}\sim_\omega\gamma_1$. Then it
  has been shown that $i_{\omega}(\gamma_{0}^{m})=i_{\omega}(\gamma_{1}^{m})$ for any $\omega\in S^{1}=\{z\in\C\mid|z|=1\}.$ Based on this fact, we always
  assume without loss of generality that each $P_{c}$ appearing in the sequel has the form (\ref{nmf}).

 \begin{lemma}\label{orient} {\rm (cf. Theorem 8.3.1 and Chap. 12 of \cite{lo2002})}  Let $c$ be an orientable closed geodesic on an $N$-dimensional Finsler manifold with its Poincar$\acute{e}$ map $P_c$.  Then, there
exists a continuous symplecitic path  $\ga$ with $\ga(0)=I$ and $\ga(1)=P_c$ such that
\bea {\rm ind}(c^{m})=i(\ga^m)
&=& m(i(\ga)+p_-+p_0-r )  -(p_- + p_0+r) - {{1+(-1)^m}\over 2}(q_0+q_+) \nn\\&&
              + 2\sum_{j=1}^r{E}\left(\frac{m\th_j}{2\pi}\right)+ 2\sum_{j=1}^{r_{\ast}}\vf\left(\frac{m\aa_j}{2\pi}\right) - 2r_{\ast},
\lb{preind}\eea
and
\bea
{\rm null}(c^{m})=\nu(\ga^m)
 &=& \nu(\ga) + {{1+(-1)^m}\over 2}(q_-+2q_0+q_+) + 2\vs(c,m),    \lb{pnul}\eea
where we denote by
\be \vs(c,m) = \left(r - \sum_{j=1}^r\vf\left(\frac{m\th_j}{2\pi}\right)\right)
             + \left(r_{\ast} - \sum_{j=1}^{r_{\ast}}\vf\left(\frac{m\aa_j}{2\pi}\right)\right)
             + \left(r_0 - \sum_{j=1}^{r_0}\vf\left(\frac{m\bb_j}{2\pi}\right)\right).    \nn\lb{yuxiang}\ee
\end{lemma}

 From now on, we focus on the non-orientable case.
\begin{lemma}\label{unorient} Let $c$ be a non-orientable closed geodesic on a $d$-dimensional Finsler manifold with its linear
Poincar$\acute{e}$ map $P_c$.  Then, the following two claims hold.

$(i)$ If $d$ is even, there is a symplectic path $\gamma$ in Sp$(2d-2)$  with $\gamma(0)=I$ and $\gamma(1)=P_c$ satisfying
$$  ({\rm ind}(c^m),{\rm null}(c^{m}))= \left\{\begin{array}{ll}
    (i_{-1}(\gamma^m),\nu_{-1}(\gamma^{m})),\; &\ if\ m\ is\ odd,\\
    (i(\gamma^m),\nu(\gamma^{m}))\; &\ if \ m  \ is \ even.\\
    \end{array}\right. $$

$(ii)$ If $d$ is odd, there is a symplectic path $\td\gamma$ in Sp$(2d)$ with $\td\gamma(0)=I$ and $\td\gamma(1)=N_{1}(1,1)\diamond P_c$ satisfying
$$
    ({\rm ind}(c^m),{\rm null}(c^{m}))= \left\{\begin{array}{ll}
    (i_{-1}(\td{\gamma}^m),\nu_{-1}(\td{\gamma}^{m})-1),\; &\ if\ m\ is\ odd,\\
    (i(\td{\gamma}^m),\nu(\td{\gamma}^{m})-1)\; &\ if \ m  \ is \ even.\\
    \end{array}\right. $$

   \end{lemma}

\Proof For the case of $m=1$, such a conclusion has been obtained by Theorem 1.1 of \cite{Liu2005}. Based on it,  Lemma \ref{unorient} is a direct application of Bott formulae (cf. Theorem 9.2.1 in \cite{lo2002}). $\hfill\Box$

For any $m\in\N$, we define
$$E_{m}(a)=E\left(a-{1-(-1)^{m}\over4}\right),\ \varphi_{m}(a)=\varphi\left(a-{1-(-1)^{m}\over4}\right),\ \forall a\in\R.$$

By Lemmas \ref{orient} and \ref{unorient}, we now derive the precise iteration formulae of Morse indices
for a non-orientable closed geodesic on a Finsler manifold .

\begin{theorem}
\label{indnullity}
 Let $c$ be a non-orientable closed geodesic on a $d$-dimensional Finsler Manifold $M$ with its linear
Poincar$\acute{e}$ map $P_c$. Then for every $m\in\N$, we have
\bea
{\rm ind}(c^{m})&=&m({\rm ind}(c)+q_{0}+q_{+}-2r')-(q_{0}+q_{+})-{1+(-1)^{m}\over2}\left(r+p_{-}+p_{0}+{1-(-1)^{d}\over2}\right)\nn\\
& &+2\sum_{j=1}^{r} E_{m}\left({m\theta_{j}\over2\pi}\right)+2\sum_{j=1}^{r_*} \varphi_{m}\left({m\alpha_{j}\over2\pi}\right)-2r_*,
\lb{ind1}\eea
and
\bea
{\rm null}(c^{m})&=&{\rm null}(c)+{1+(-1)^{m}\over2}\left(p_{-}+2p_{0}+p_{+}+{1-(-1)^{d}\over2}\right)+2\td{\varsigma}(c,m),\lb{null1}
\eea
where we denote by
$$\td{\varsigma}(c,m)=\left(r-\sum_{j=1}^{r}\varphi_{m}\left({m\theta_{j}\over2\pi}\right)\right)
+\left(r_{*}-\sum_{j=1}^{r_{*}}\varphi_{m}\left({m\alpha_{j}\over2\pi}\right)\right)
+\left(r_{0}-\sum_{j=1}^{r_{0}}\varphi_{m}\left({m\beta_{j}\over2\pi}\right)\right).$$
\end{theorem}
 \Proof We only prove the case when $d$ is even and $m$ is odd, since it is just the case we encounter in this paper and the other cases can be proved similarly. By Lemma \ref{unorient}, there exists a symplectic path $\gamma$ in Sp$(2d-2)$ with $\ga(0)=I$ and $\gamma(1)=P_c$ such that
\bea (\text{ind}(c^{m}),\text{null}(c^{m}))=(i_{-1}(\gamma^{m}),\nu_{-1}(\gamma^{m})), \forall\ m\in2\N-1.
\lb{201679a}\eea
It together with the Bott-type formulae (cf. Theorem 9.2.1 of \cite{lo2002}) and Lemma \ref{orient} gives
\bea i_{-1}(\gamma^{m})&=&i(\gamma^{2m})-i(\gamma^{m})\nn\\
&=&m(i(\ga)+p_-+p_0-r )-(q_0+q_+)    \nn\\
&&+ 2\sum_{j=1}^r\left[{E}\left(\frac{m\th_j}{\pi}\right)-{E}\left(\frac{m\th_j}{2\pi}\right)\right]+ 2\sum_{j=1}^{r_{\ast}}\left[\vf\left(\frac{m\aa_j}{\pi}\right)-\vf\left(\frac{m\aa_j}{2\pi}\right)\right]\nn\\
&=&m(i(\ga)+p_-+p_0-r )-(q_0+q_+)  \nn\\
&&+ 2\sum_{j=1}^r{E}\left(\frac{m\th_j}{2\pi}-{1\over2}\right)+ 2\sum_{j=1}^{r_{\ast}} \vf\left(\frac{m\aa_j}{2\pi}-{1\over2}\right)-2r_{*}\lb{201679b}\\
&=&m(i_{-1}(\gamma)+q_0+q_+-2r' )-(q_0+q_+) \nn\\
&&+ 2\sum_{j=1}^r{E}\left(\frac{m\th_j}{2\pi}-{1\over2}\right)+ 2\sum_{j=1}^{r_{\ast}} \vf\left(\frac{m\aa_j}{2\pi}-{1\over2}\right)-2r_{*},\nn
\eea
where the third identity we have used $E(2a)-E(a)=E\left(a-{1\over2}\right)$ and $\vf(2a)-\vf(a)=\vf\left(a-{1\over2}\right)-1,$ and the last identity is due to
$$i(\gamma)=i_{-1}(\gamma)+(q_{0}+q_{+})+(r-2r')-(p_{0}+p_{-}),$$ which is a result of direct computation on splitting numbers based on Theorem 12.2.3 of \cite{lo2002}.

Observing by definition $\nu_{-1}(\gamma)=q_{-}+2q_{0}+q_{+}$, we obtain similarly as above that
 \bea \nu_{-1}(\gamma^{m})&=&\nu(\gamma^{2m})-\nu(\gamma^{m})\nn\\
&=&(q_{-}+2q_{0}+q_{+})-2\sum_{j=1}^r\left[{\varphi}\left(\frac{m\th_j}{\pi}\right)-{\varphi}\left(\frac{m\th_j}{2\pi}\right)\right]\nn\\
&&-2\sum_{j=1}^{r_{\ast}}\left[\vf\left(\frac{m\aa_j}{\pi}\right)-\vf\left(\frac{m\aa_j}{2\pi}\right)\right]
-2\sum_{j=1}^{r_{0}}\left[\vf\left(\frac{m\bb_j}{\pi}\right)-\vf\left(\frac{m\bb_j}{2\pi}\right)\right]\nn\\
&=&\nu_{-1}(\gamma)+2\left(r-\sum_{j=1}^r{\varphi}\left(\frac{m\th_j}{2\pi}-{1\over2}\right)\right)\lb{201679c} \\
&& +2\left(r_{*}-\sum_{j=1}^{r_{\ast}}\vf\left(\frac{m\aa_j}{2\pi}-{1\over2}\right)\right)
+2\left(r_0-\sum_{j=1}^{r_{0}}\vf\left(\frac{m\bb_j}{2\pi}-{1\over2}\right)\right).\nn
\eea

Thus (\ref{ind1}) and (\ref{null1})  immediately follow from (\ref{201679a}), (\ref{201679b}) and (\ref{201679c}). \hfill$\Box$

\subsection{A variant of Precise index iteration formulae}

In this section, we give a variant of the precise index iteration formulae in section 3.1 which makes them more intuitive and enables us to apply
the Kronecker's approximation theorem to study the multiplicity of non-contractible closed geodesics on $\R P^{n}$.

To prove Theorem \ref{mainresult}, we always assume that there exists only one non-contractible prime closed geodesic
$c$ on $M=\RP$ with a bumpy metric $F$, which is then just the well known minimal point of the energy functional
$E$ on $\Lambda_{g}M$ satisfying ${\rm ind}(c)=0$. Now the Morse-type number is given by
\begin{equation*} m_q \equiv M_q(\Lm_gM)=\sum_{m\ge 1}\dim{\ol{C}}_q(E, c^{2m-1}), \quad
                \forall q\in \N\cup\{0\}. \nn\end{equation*}
Then by Lemma \ref{Rad1992}(i), Lemma \ref{Lm2.3} and Morse inequality, we have the following conclusion.
\begin{lemma}
\label{morsebetti} {\rm (cf. Lemma 3.1 of \cite{DLX2015})} Assuming the existence of only one non-contractible
prime closed geodesic $c$ on $\RP$ with a bumpy metric $F$, there hold
\begin{equation} m_{2q+1}=\bar{\beta}_{2q+1}=0\  {\rm and}\ m_{2q}=\bar{\beta}_{2q}, \qquad
           \forall\ q\in \N\cup\{0\}.{\label{beqm}}\end{equation}
\end{lemma}

We consider two cases according to the parity of dimension of the real projective space.
First we study the case of $\R P^{2n+1}$. Note that the other one behaves similarly.

\begin{lemma}\label{lemma3.2} Suppose $c$ is the only one non-contractible
prime closed geodesic $c$ on $(\R P^{2n+1}, F)$ with a bumpy metric $F$. Then there exist
$\hat{\theta}_1$, $\hat{\theta}_2$, \dots, $\hat{\theta}_k$ in $\Q^c$ with $2\leq k\leq2n$ such that
\bea  \sum_{j=1}^{k}\hat{\theta}_{j}&=&{1\over2}\left(k+{n\over n+1}\right), \label{mean}\\
 {\rm ind}(c^{m})&=& m\left({n\over n+1}\right)+k-2\sum_{j=1}^{k}\left\{{m\hat{\theta}_{j}}\right\},\qquad\forall\ m\ge 1. \lb{iter1}\eea
\end{lemma}

\Proof See (3.6), (3.7) and (3.8) in \cite{DLX2015}. Also compare the proof of Lemma 3.6.\hfill$\Box$

Now we give a variant of the precise index iteration formulae (\ref{iter1}) specially for our purpose. Let  $m=2(n+1)l+2L+1$
with $l\in\N$ and $L\in\Z$. By (\ref{mean}) and (\ref{iter1}) we obtain
\bea  \text{ind}(c^{m})
&=& 2nl+k+(2L+1){n\over n+1 }\nn\\
& & -2\left(\left\{{k\over2}+{(2L+1)n\over2(n+1)}-\sum_{j=2}^{k}\left\{m\hat{\theta}_{j}\right\}\right\}
    +\sum_{j=2}^{k}\left\{m\hat{\theta_{j}}\right\}\right)  \nn\\
&=& 2nl+2\left[{k\over2}+{(2L+1)n\over2(n+1)}\right]+2\left\{{k\over2}+{(2L+1)n\over2(n+1)}\right\}  \nn\\
& & -2\left(\left\{\left\{{k\over2}+{(2L+1)n\over2(n+1)}\right\}-\sum_{j=2}^{k}\left\{m\hat{\theta}_{j}\right\}\right\}
   +\sum_{j=2}^{k}\left\{m\hat{\theta_{j}}\right\}\right)  \nn\\
&=& 2nl+2\left[Q_{L}\right]+2\left\{Q_{L}\right\}-2\left(\left\{\left\{Q_{L}\right\}-\sum_{j=2}^{k}
  \left\{m\hat{\theta}_{j}\right\}\right\}+\sum_{j=2}^{k}\left\{m\hat{\theta_{j}}\right\}\right),
  \label{evenformulae}\eea
where in the last identity for notational simplicity, we denote by $Q_{L}={k\over2}+{(2L+1)n\over2(n+1)}. $

Since $\sum_{j=2}^{k}\{m\hat{\theta}_{j}\}\in\Q^{c}$, we obtain by (\ref{evenformulae}) that for $1\leq i\leq k-2$,
\be  \text{ind}(c^{m})
= \left\{\begin{array}{ll}
2nl+2\left[Q_{L}\right],& \text{iff}\ \sum_{j=2}^{k}\left\{m\hat{\theta}_{j}\right\}\in (0,\left\{Q_{L}\right\}),\\
2nl+2\left[Q_{L}\right]-2i,& \text{iff}\ \sum_{j=2}^{k}\left\{m\hat{\theta}_{j}\right\}\in (i-1+\left\{Q_{L}\right\},i+\left\{Q_{L}\right\}),\\
2nl+2\left[Q_{L}\right]-2(k-1),& \text{iff}\ \sum_{j=2}^{k}\left\{m\hat{\theta}_{j}\right\}\in (k-2+\left\{Q_{L}\right\},k-1).
\end{array}\right.   \lb{evenodddengjia}\ee
Let $I_{0}(L)=(0,\left\{Q_{L}\right\})$, $I_{k-1}(L)=(k-2+\left\{Q_{L}\right\},k-1),$ and
$$  I_{i}(L)=(i-1+\left\{Q_{L}\right\},i+\left\{Q_{L}\right\})\quad \text{for}\quad 1\le i\le k-2.  $$
Then, (\ref{evenodddengjia}) can be stated in short as that for any integers $m=2(n+1)l+2L+1$ and  $0\leq i\leq k-1,$
\be   \text{ind}(c^{m}) = 2nl+2\left[Q_{L}\right]-2i\quad \text{if and only if}\quad
    \sum_{j=2}^{k}\left\{m\hat{\theta}_{j}\right\}\in I_{i}(L).   \lb{shortdengjia}\ee

\begin{remark}
\label{transformation} Let $(\tau(1),\tau(2),\dots,\tau(k))$ be an arbitrary permutation of $(1,2,\dots,k)$. Then, the same conclusion as (\ref{shortdengjia}) with $j$ ranging in $\{\tau(1),\tau(2),\dots,\tau(k-1)\}$ instead is still valid.
\end{remark}

The following lemma will be also needed in the proof of Theorem \ref{mainresult} for $\R P^{2n+1}$  in Section 5.
\begin{lemma}\label{bound} Under the assumption of Lemma \ref{lemma3.2}, for any positive integers $l$ and $m$, we have
$$  |{\rm ind}(c^{m})-2nl|>2n\quad \text{holds whenever}\quad  |m-2(n+1)l|>4(n+1). $$
\end{lemma}

\Proof From (\ref{iter1}), we have
$$  {\rm ind}(c^{m}) = 2nl+(m-2(n+1)l)\cdot{n\over n+1}+k-2\sum_{j=1}^{k}\left\{{m\hat{\theta}_{j}}\right\}, $$
which yields immediately that
\bea  |{\rm ind}(c^{m})-2nl|
&\ge& |m-2(n+1)l|\cdot{n\over n+1}-|k-2\sum_{j=1}^{k}\left\{{m\hat{\theta}_{j}}\right\}|  \nn\\
&>& 4n-k \ge 4n-2n = 2n,  \nn\eea
where the fact $k\le 2n$ is used. \hfill$\Box$

For the case of $\R P^{2n}$, similar to Lemma \ref{lemma3.2}, we have

\begin{lemma}\label{lemma3.4} Suppose $c$ is the only one non-contractible
prime closed geodesic $c$ on $(\R P^{2n}, F)$ with a bumpy metric $F$. Then, there exist
$\hat{\theta}_1$, $\hat{\theta}_2$, \dots, $\hat{\theta}_{2r}$ in $\Q^c$ with $2\leq r\leq2n-1$ such that
\bea  \sum_{j=1}^{2r}\hat{\theta}_{j}&=&{1\over2}\left(2r+{2n-1\over 2n}\right), \label{mean33}\\
 {\rm ind}(c^{2m-1})&=& (2m-1)\left({2n-1\over 2n}\right)+2r-2\sum_{j=1}^{2r}\left\{{(2m-1)\hat{\theta}_{j}}\right\},\qquad\forall\ m\ge 1. \lb{iter33}\eea
\end{lemma}

\Proof Since the Finsler metric F is bumpy, it follows $\text{null}(c^{m})=0$ for every $m\in\N$. In particular, $\text{null}(c)=\nu_{-1}(\gamma)=q_{-}+2q_{0}+q_{+}=0$, which implies $q_-=q_0=q_+=0.$ In addition by (\ref{null1}), $\text{null}(c^{2})=0$ then yields $p_-=p_0=p_+=0.$ As a result, we get
$$\td{\varsigma}(c,m)=0,\ \forall m\in\N,$$
and so $\frac{\th_{j}}{2\pi}'s$, $\frac{\aa_{j}}{2\pi}'s$ and $\frac{\bb_{j}}{2\pi}'s$ are all in $\Q^c\cap(0,1)$.

Due to ${\rm ind}(c)=0$, by (\ref{ind1}) in Theorem \ref{indnullity} we obtain
\bea {\rm ind}(c^{2m-1})
&=&-2(2m-1)r'+\sum_{j=1}^r E\left({(2m-1)\theta_{j}\over2\pi}-{1\over2}\right)\nn\\
&=&-2(2m-1)r'+2\sum_{j=1}^r \left(\left({(2m-1)\theta_{j}\over2\pi}-{1\over2}\right)-\left\{{(2m-1)\theta_{j}\over2\pi}-{1\over2}\right\}+1\right)\nn\\
&=&(2m-1)\left(-2r'+\sum_{j=1}^{r}{\theta_{j}\over\pi}\right)-2\sum_{j=1}^r \left(\left\{{(2m-1)\theta_{j}\over2\pi}-{1\over2}\right\}-{1\over2}\right)\nn\\
&=&(2m-1)\left(-2r'+\sum_{j=1}^{r}{\theta_{j}\over\pi}\right)+2r-2\sum_{j=1}^r \left(\left\{{(2m-1)\theta_{j}\over\pi}\right\}+\left\{{-(2m-1)\theta_{j}\over2\pi}\right\}\right)\nn\\
\label{iter44}\eea
which implies
$ \hat{i}(c)=-2r'+\sum_{j=1}^{r}{\theta_{j}\over\pi}.$
It together with (\ref{breident1}) of Theorem \ref{Thm1.1} yields
\bea \sum_{j=1}^r\frac{\th_j}{2\pi}
=\frac{1}{2}\left(2r'+\frac{2n-1}{2n}\right).\label{mean44}\eea

Let $\hat{\theta}_j=\frac{\th_j}{\pi}-[\frac{\th_j}{\pi}]+1$ for $1\leq j \leq r$ and $\hat{\theta}_j=-\frac{\th_{j-r}}{2\pi}$ for $r+1\leq j \leq 2r$,
then (\ref{mean33}) follows from (\ref{mean44}), and (\ref{iter33}) follows from (\ref{iter44}) and (\ref{mean44}).
\hfill$\Box$

\begin{remark}
If we replace $2n-1$ and $2r$ in Lemma \ref{lemma3.4} by $n$ and $k$ respectively, (\ref{mean33}) and (\ref{iter33}) are just the same
form as (\ref{mean}) and (\ref{iter1}) respectively. Hence (\ref{evenformulae})-(\ref{shortdengjia}) also hold when we replace
$n$ and $k$ by $2n-1$ and $2r$ respectively.
\end{remark}

For the case of $\R P^{2n}$, similar to Lemma \ref{bound}, we have

\begin{lemma}\label{bound2} Under the assumption of Lemma \ref{lemma3.4}, for any positive integers $l$ and $m\in2\N-1$, we have
$$  |{\rm ind}(c^{m})-2(2n-1)l|>2(2n-1)\quad \text{holds whenever}\quad  |m-4nl|>8n. $$
\end{lemma}


\section{The system of irrational numbers}

Let $\aa=\{\alpha_{1}, \alpha_{2}, \ldots, \alpha_{m}\}$ be a set of $m$ irrational numbers. As usual,
we have
\begin{definition}\lb{rank}
The set $\aa$ of irrational numbers is linearly independent over $\Q$, if there do not exist $c_1$,
$c_2$, $\ldots$, $c_m$ in $\Q$ such that $\sum_{j=1}^{m}|c_{j}|> 0$ and
\be  \sum_{j=1}^{m}c_{j}\alpha_{j}\in\Q,  \label{81a}\ee
and is linearly dependent over $\Q$ otherwise. The rank of $\aa$ is defined to be the number of
elements in a maximal linearly independent subset of $\aa$, which we denote by $\rank(\aa)$.
\end{definition}

\begin{lemma}\lb{L3.1}
Let $r=\rank(\aa)$. Then there exist $p_{jl}\in\Z$, $\beta_{l}\in\Q^{c}$ and $\xi_{j}\in \Q$ for $1\leq j\leq r$ and $1\leq l\leq m$ such that
 \begin{equation}
\label{81b}
\alpha_{j}=\sum_{l=1}^{r}p_{jl}\beta_{l}+\xi_{j},\ \forall 1\leq j\leq m.
\end{equation}

\end{lemma}

\Proof  Let $\aa'=\{\aa_{m_1}, \aa_{m_2}, \ldots, \aa_{m_r}\}$ be a maximal linearly independent subset
of $\aa$. Then there exist $c_{jl}\in\Q$ and $\xi_{j}\in \Q$ such that
 \begin{equation}
\label{81c}
\alpha_{j}=\sum_{l=1}^{r}c_{jl}\alpha_{ml}+\xi_{j},\ \forall 1\leq j\leq m.
\end{equation}
 For every $1\leq l\leq r$, we define $J_{l}=\{1\leq j\leq m\mid c_{jl}\neq 0\}$ and then for $j\in J_{l}$ let $c_{jl}={r_{jl}\over q_{jl}}$ with $r_{jl}$ prime to $q_{jl}$. Define
 $q_{l}=\prod_{j\in J_{l}} q_{jl}$ and $$\beta_{l}={\alpha_{ml}\over q_{l}}\in\Q^{c}\ \text{and}\ p_{jl}=q_{l}c_{jl}\in\Z,\ \forall 1\leq j\leq m\ \text{and}\ 1\leq l\leq r.$$
 Then, (\ref{81b}) follows. \hfill$\Box$

In order to study the multiplicity of closed geodesics on $(\R P^{2n+1},F)$ with a bumpy Finsler metric
$F$, we are particularly interested in the irrational system  $\{\hat{\theta}_{1}$,
$\hat{\theta}_{2}$,\dots,$\hat{\theta}_{k}\}$ with rank $1$ satisfying (\ref{mean}).
Then by Lemma \ref{L3.1}, it can be reduced to the following system
\be   \hat{\theta}_{j}=p_{j}\theta+\xi_{j},\ \forall 1\leq j\leq k,  \label{Qdependent}\ee
with $\theta\in\mathbb{Q}^{c}$, $p_{j}\in \Z\backslash\{0\}$, $\xi_{j}\in\mathbb{Q}\cap[0,1)$
satisfying
\bea
&& p_{1}+p_{2}+\cdots+p_{k}=0, \lb{basic1c}\\
&& \{\xi_{1}+\xi_{2}+\cdots+\xi_{k}\}\in(0,1)\backslash\{{1/2}\}, \lb{basic2c}\eea
where to get $\xi_j\in [0,1)$, if necessary, we can replace $\hat{\th}_j$ and $\xi_j$ by
$\breve{\theta}_{j}=\hat{\theta}_{j}-[\xi_{j}]$ and $\breve{\xi}_{j}=\{\xi_{j}\}.$

Take arbitrarily $\eta\in\Q$ and make  the following natural $\eta$-action  to the system (\ref{Qdependent}):
\begin{equation}
\label{eta}
{\eta}(\theta)=\theta+\eta,\ {\eta}(\hat{\theta}_{j})=\hat{\theta}_{j}-\left[\xi_{j}-p_{j}\eta\right]\ \text{and}\ {\eta}({\xi}_{j})=\left\{\xi_{j}-p_{j}\eta\right\},\ \forall 1\leq j\leq k,
\end{equation}
which is obviously induced by the transformation ${\eta}(\theta)=\theta+\eta$.
Then, we get a new system
\begin{equation}
\label{Qdependent1}
{\eta}(\hat{\theta}_{j})=p_{j}{\eta}(\theta)+{\eta}(\xi_{j}),\ \forall 1\leq j\leq k,
\end{equation}
with
\begin{equation}
\label{invariant1}
\begin{aligned}
\{{\eta}({\xi}_{1})+{\eta}({\xi}_{2})+\cdots+{\eta}({\xi}_{k})\}&=\{\{\xi_{1}-p_{1}\eta\}+\{\xi_{2}-p_{2}\eta\}+\cdots+\{\xi_{k}-p_{k}\eta\}\}\\
&=\{\xi_{1}+\xi_{2}+\cdots+\xi_{k}-(p_{1}+p_{2}+\dots+p_{k})\eta\}\\
&=\{\xi_{1}+\xi_{2}+\cdots+\xi_{k}\},
\end{aligned}
\end{equation}
where the third equality we have used the condition (\ref{basic1c}).
For simplicity of writing, we also denote  the new system (\ref{Qdependent1})  by (\ref{Qdependent})$_{\eta}$ meaning that it comes from
(\ref{Qdependent}) by an $\eta$-action.

For the system (\ref{Qdependent})$_{\eta}$ with $\eta\in\Q$, we divide the set $\{1\leq j\leq k\}$ into the following three parts:
\bea \mathcal{K}_{0}^{+}(\eta)&=&\{1\leq j\leq k\mid \eta(\xi_{j})=0,\ p_{j}>0\}, \nn\\
\mathcal{K}_{0}^{-}(\eta)&=&\{1\leq j\leq k\mid \eta(\xi_{j})=0,\ p_{j}<0\}, \nn\\
\mathcal{K}_{1}(\eta)&=&\{1\leq j\leq k\mid \eta(\xi_{j})\neq0 \}.\eea
Denote by $k_{0}^{+}(\eta)$, $k_{0}^{-}(\eta)$ and $k_{1}(\eta)$ the numbers ${}^{\#}\mathcal{K}_{0}^{+}(\eta)$, ${}^{\#}\mathcal{K}_{0}^{-}(\eta)$ and ${}^{\#}\mathcal{K}_{1}(\eta)$ respectively. For the case of $\eta=0$, we write them for short as $k_{0}^{+}$, $k_{0}^{-}$ and $k_{1}.$ It follows immediately that $$k_{0}^{+}(\eta)+k_{0}^{-}(\eta)+k_{1}(\eta)=k.$$
By (\ref{basic2c}) and (\ref{invariant1}), it is obvious that $k_{1}(\eta)\geq1$ for every $\eta\in\Q$.
\begin{definition}
\label{crucialnumber}
For every $\eta\in\Q$, the {\bf absolute difference number} of (\ref{Qdependent})$_{\eta}$ is defined to be the non-negative number $|k_{0}^{+}(\eta)-k_{0}^{-}(\eta)|.$
 The {\bf effective difference number} of {(\ref{Qdependent})} is defined by
$$\max\{|k_{0}^{+}(\eta)-k_{0}^{-}(\eta)|\mid \eta\in\Q\}.$$
Two  systems of irrational numbers with rank $1$ are called to be {\bf equivalent}, if their effective difference numbers are the same one.
\end{definition}

\begin{remark} By the definition of an $\eta$-action in (\ref{eta}), it can be checked directly that
 $\eta_{1}\circ\eta_{2}=\eta_{1}+\eta_{2}$  for every $\eta_{1}$ and $\eta_{2}$ in $\Q.$ So every system of irrational numbers with rank $1$ is equivalent to the one which comes from itself  by an $\eta$-action.
\end{remark}

We have first the following simple equivalent pairs.
\begin{lemma}
\label{equivalentsystem}
Assume that
\begin{equation}
\label{20150729a}
\left\{\begin{array}{ll}
\hat{\theta}_{j}=p_{j}\theta+\xi_{j},& \forall 1\leq j\leq k-1,\\
\hat{\theta}_{k}=p_{k}\theta,\\
\end{array}
\right.
\end{equation}
with $\sum_{j=1}^{k}p_{k}=0$ and $\left\{\sum_{j=1}^{k-1}\xi_{k}\right\}\in(0,1)\backslash\{1/2\}.$

Then, (\ref{20150729a}) is equivalent to
\begin{equation}
\label{20150729b}
\left\{\begin{array}{ll}
\hat{\theta}_{j}=p_{j}\theta+\xi_{j},& \forall\ 1\leq j\leq k-1,\\
\hat{\theta}_{k,l}=\text{sgn}(p_{k})\theta+{l\over |p_{k}|},& \forall\ 0\leq l\leq |p_{k}|-1,\\
\end{array}
\right.
\end{equation}
 where as usual we define $\sgn(a)=\pm 1$ for $a\in\R\bs\{0\}$ when $\pm a>0$.
\end{lemma}
\Proof Take $\eta\in\Q$ arbitrarily and recall the definition of $\eta$-action in (\ref{eta}). Then the equation $\hat{\theta}_{k}=p_{k}\theta$ contributes $\text{sgn}(p_{k})$ to the absolute difference number of (\ref{20150729a})$_{\eta}$ if and only if
$$\eta(0)=\{0-p_{k}\eta\}=\{-p_{k}\eta\}=0,$$ that is $\eta\in\Z_{|p_{k}|}$, which is also the sufficient and necessary condition such that the equations
 $$\hat{\theta}_{k,l}=\text{sgn}(p_{k})\theta+{l\over |p_{k}|},\ \forall\ 0\leq l\leq |p_{k}|-1,$$
 contribute $\text{sgn}(p_{k})$ to the absolute difference number of (\ref{20150729b})$_{\eta}$. Since the other equations with $1\leq j\leq k-1$  in (\ref{20150729a}) and (\ref{20150729b}) are the same, so do their contributions to the absolute difference numbers of (\ref{20150729a})$_{\eta}$ and (\ref{20150729b})$_{\eta}$.
 As a result, the absolute difference numbers of (\ref{20150729a})$_{\eta}$ and (\ref{20150729b})$_{\eta}$ are equal for any $\eta\in\Q$ which yields that the effective difference numbers of (\ref{20150729a}) and (\ref{20150729b}) are the same and so they are equivalent. \hfill$\Box$
\begin{remark}
\label{beautifulchange}
For the system (\ref{20150729b}), we have
$$\left\{\sum_{j=1}^{k-1}\xi_{j}+\sum_{l=0}^{|p_{k}|-1}{l\over |p_{k}|}\right\}
=\left\{\begin{array}{ll}
\left\{\sum_{j=1}^{k-1}\xi_{j}\right\},& if\ p_{k}\ is\ odd,\\
\left\{\sum_{j=1}^{k-1}\xi_{j}+{1\over2}\right\},& if\ p_{k}\ is\ even.
\end{array}
\right.$$
By the assumption of $\left\{\sum_{j=1}^{k-1}\xi_{j}\right\}\in(0,1)\backslash\{1/2\}$, it follows that
$$\left\{\sum_{j=1}^{k-1}\xi_{j}+\sum_{l=0}^{|p_{k}|-1}{l\over |p_{k}|}\right\}\in(0,1)\backslash\{1/2\}.$$
\end{remark}

\begin{lemma}
\label{cutoff}
If there exist $1\leq j'<j''\leq k$ satisfying that $p_{j'}\cdot p_{j''}=-1$ and
$\left\{{\xi_{j'}}+ {\xi_{j''}}\right\}=0$ in
\begin{equation}
\label{Qdependent1029}
\hat{\theta}_{j}=p_{j}\theta+\xi_{j},\ \forall 1\leq j\leq k,
 \end{equation}
 then (\ref{Qdependent1029}) is equivalent to the system
\begin{equation}
\label{cutoffequations}
\hat{\theta}_{j}=p_{j}\theta+\xi_{j},\ \forall j\in\{1,2,\dots,k\}\backslash\{j',j''\}.
\end{equation}
\end{lemma}
\Proof Assume without loss of generality that $p_{j'}=-p_{j''}=1$ and take $\eta\in\Q$  arbitrarily. Then by (\ref{eta}) and the given condition, we have
$$
\left\{{\eta}(\xi_{j'})+{\eta}(\xi_{j''})\right\}=\left\{ \{\xi_{j'}-\eta\} + \{\xi_{j''}+\eta\} \right\}=\{\xi_{j'}+\xi_{j''}\}=0.
$$
Thus, ${\eta}({\xi_{j'}})=0$ if and only if ${\eta}({\xi_{j''}})=0,$ that is, $j'\in\mathcal{K}_{0}^{+}({\eta})$ if and only if
$j''\in\mathcal{K}_{0}^{-}({\eta}).$ As a result, $p_{j'}$ and $p_{j''}$ together contribute nothing to the absolute difference number
of (\ref{Qdependent1029})$_\eta$ for any $\eta\in\Q$. It then follows immediately that (\ref{Qdependent1029}) is equivalent to
(\ref{cutoffequations}). \hfill$\Box$

The following theorem is our main result of this section which is concerned with the lower estimate on the effective difference number
of (\ref{Qdependent}) and will play a crucial role in our proof of Theorem \ref{mainresult} in Section 5.
\begin{theorem}
\label{crucialtheorem}
For every system of irrational numbers (\ref{Qdependent}) satisfying the conditions (\ref{basic1c}) and (\ref{basic2c}),  it holds that
   \begin{equation}
   \label{crucialestimate}
   \max\{|k_{0}^{+}(\eta)-k_{0}^{-}(\eta)|\mid \eta\in\Q\}\geq1.
   \end{equation}
\end{theorem}

\begin{remark}
The condition (\ref{basic2c}) can not be replaced by the  weaker condition
\begin{equation}
\label{weakcondition}
\{\xi_{1}+\xi_{2}+\cdots+\xi_{k}\}\in(0,1).
\end{equation} For instance, we consider the system
$\hat{\theta}_{1}=-\theta+{1\over2},\ \hat{\theta}_{2}=-\theta,\ \hat{\theta}_{3}=2\theta,$
which satisfies  the conditions (\ref{basic1c}) and (\ref{weakcondition}) but (\ref{basic2c}).
However, one can check directly that $|k_{0}^{+}(\eta)-k_{0}^{-}(\eta)|=0$ for any $\eta\in\Q.$
As we will see, such a phenomenon does not occur if  the condition (\ref{basic2c}) holds.
\end{remark}

{\bf Proof of Theorem \ref{crucialtheorem}:} We carry out the proof with two steps.

{\bf Step 1:} First, letting $\eta_{k}={\xi_{k}\over p_{k}}$ and making ${\eta_{k}}$-action to the original system (\ref{Qdependent}), we obtain by (\ref{eta}) that
\begin{equation}
\label{729change1}
\left\{\begin{array}{ll}
{\eta_{k}}(\hat{\theta}_{j})=p_{j}{\eta_{k}}(\theta)+{\eta_{k}}(\xi_{j}),\ \forall 1\leq j\leq k-1,\\
{\eta_{k}}(\hat{\theta}_{k})=p_{k}{\eta_{k}}(\theta).
\end{array}
\right.
\end{equation}
Then by Lemma \ref{equivalentsystem}, the system (\ref{729change1}) is equivalent to
\begin{equation}
\label{729change11}
\left\{\begin{array}{ll}
{\eta_{k}}(\hat{\theta}_{j})=p_{j}{\eta_{k}}(\theta)+{\eta_{k}}(\xi_{j}),& \forall\ 1\leq j\leq k-1,\\
\hat{\theta}_{k,l'}=\text{sgn}(p_{k}){\eta_{k}}(\theta)+{l'\over |p_{k}|},& \forall\ 0\leq l'\leq |p_{k}|-1,\\
\end{array}
\right.
\end{equation}

Secondly, taking $\eta_{k-1}\in\Q$ such that ${\eta_{k-1}}\circ{\eta_{k}}(\xi_{k-1})=0$ and making ${\eta_{k-1}}$-action to the system (\ref{729change11}), we get
\begin{equation}
\label{729change2}
\left\{\begin{array}{ll}
{\eta_{k-1}}\circ{\eta_{k}}(\hat{\theta}_{j})=p_{j}{\eta_{k-1}}\circ{\eta_{k}}(\theta)+{\eta_{k-1}}\circ{\eta_{k}}(\xi_{j}),& \forall\ 1\leq j\leq k-2,\\
{\eta_{k-1}}\circ{\eta_{k}}(\hat{\theta}_{k-1})=p_{k-1}{\eta_{k-1}}\circ{\eta_{k}}(\theta),& \\
{\eta_{k-1}}\circ(\hat{\theta}_{k,l'})=\text{sgn}(p_{k}){\eta_{k-1}}\circ{\eta_{k}}(\theta)+{\eta_{k-1}}({l'\over |p_{k}|}),& \forall\ 0\leq l'\leq |p_{k}|-1.\\
\end{array}
\right.
\end{equation}
 Again by Lemma \ref{equivalentsystem}, the system (\ref{729change2}) is equivalent to
\begin{equation}
\label{729change22}
\left\{\begin{array}{ll}
{\eta_{k-1}}\circ{\eta_{k}}(\hat{\theta}_{j})=p_{j}{\eta_{k-1}}\circ{\eta_{k}}(\theta)+\eta_{k-1}\circ{\eta_{k}}(\xi_{j}),& \forall\ 1\leq j\leq k-2,\\
\hat{\theta}_{k-1,l''}=\text{sgn}(p_{k-1}){\eta_{k-1}}\circ{\eta_{k}}(\theta)+{l''\over |p_{k-1}|},& \forall\ 0\leq l''\leq |p_{k-1}|-1,\\
{\eta_{k-1}}(\hat{\theta}_{k,l'})=\text{sgn}(p_{k}){\eta_{k-1}}\circ{\eta_{k}}(\theta)+{\eta_{k-1}}({l'\over |p_{k}|}),& \forall\ 0\leq l'\leq |p_{k}|-1,\\
\end{array}
\right.
\end{equation}

Repeating the above procedure for the rest equations with $j={k-2}$, ${k-3}$, $\cdots$, $2$, $1$ one at a time in order,
we can finally get a system equivalent to the original system (\ref{Qdependent}) which can be written in a simple form such as
\begin{equation}
\label{finalsystem}
\hat{\alpha}_{jl}=\text{sgn}(p_{j})\alpha+\xi_{jl},\ \forall\  1\leq j\leq k\ \text{and}\ 0\leq l\leq |p_{j}|-1,
\end{equation}
with $\alpha\in\Q^{c}$ and $\xi_{jl}\in\Q\cap[0,1)$. Moreover,
by (\ref{invariant1}) and Remark \ref{beautifulchange} we have
\begin{equation}
\label{incredible}
\left\{\sum_{j=1}^{k}\sum_{l=0}^{|p_{j}|-1}\xi_{jl}\right\}\in(0,1)\backslash\{{1/2}\}.
\end{equation}

{\bf Step 2:} We can cut off all the superfluous equations of the system (\ref{finalsystem}), if there are such pairs as that in Lemma \ref{cutoff}.
That is, (\ref{finalsystem}) is equivalent to some a system
\begin{equation}
\label{finalsystem1}
\hat{\theta}_{i}^{\prime}=p_{i}^{\prime}\alpha+\xi_{i}^{\prime},\ \forall\  1\leq i\leq \bar{k},
\end{equation}
with $|p_{i}^{\prime}|=1$, $\sum_{i=1}^{\bar{k}}p_{i}^{\prime}=0$ and
\begin{equation}
\label{basic2cc}
\left\{\sum_{i=1}^{\bar{k}}\xi_{i}^{\prime}\right\}\in(0,1)\backslash\{1/2\}.
\end{equation}
Here notice that $\bar{k}\geq1$ is ensured by the condition (\ref{basic2cc}).

Since all the superfluous equations are cut off, it follows that ${\bar{k}}_{0}^{+}\cdot {\bar{k}}_{0}^{-}=0.$ Assume without loss of generality that
${\bar{k}}_{0}^{+}={\bar{k}}_{0}^{-}=0,$  otherwise we have nothing to do. Since $\sum_{i=1}^{\bar{k}}p_{i}^{\prime}=0$,
we get $${}^{\#}\{1\leq i\leq \bar{k}\mid p_{i}^{\prime}=1\}={}^{\#}\{1\leq i\leq \bar{k}\mid p_{i}^{\prime}=-1\}.$$
Take arbitrarily out $i_{1}\in \{1\leq i\leq \bar{k}\mid p_{i}^{\prime}=1\}.$ Let $\bar{\eta}=\xi_{i_{1}}'$ and make the $\bar{\eta}$-action  to (\ref{finalsystem1}).
Then it follows immediately that  $\bar{k}_{0}^{+}(\bar{\eta})\geq1$. Recalling again that all the superfluous equations have been cut off at the beginning of Step 2,
we obtain $\bar{\eta}(\xi_{i}')=\{\xi_{i_{1}}'+\xi_{i}'\}\neq0$ for every $i\in\{1\leq i\leq \bar{k}\mid p_{i}^{\prime}=-1\}$ which yields $\bar{k}_{0}^{-}(\bar{\eta})=0.$  As a result, we get
$$\max\{|\bar{k}_{0}^{+}({\eta})-\bar{k}_{0}^{-}({\eta})\mid\eta\in\Q\}\geq|\bar{k}_{0}^{+}(\bar{\eta})-\bar{k}_{0}^{-}(\bar{\eta})|=\bar{k}_{0}^{+}(\bar{\eta})\geq1.$$
Since the original system (\ref{Qdependent}) is equivalent to (\ref{finalsystem1}), the estimate (\ref{crucialestimate}) follows immediately. \hfill$\Box$

The proof of Theorem \ref{crucialtheorem} can be illuminated by the concrete example below.

\begin{example}
Consider the irrational system
\begin{equation}
\label{simpleexample1}
\hat{\theta}_{1}=-\theta+{5\over6},\ \hat{\theta}_{2}=-2\theta+{1\over3},\ \hat{\theta}_{3}=3\theta+{1\over2}.
\end{equation}
\end{example}
One can check directly that the system (\ref{simpleexample1}) is a special case of  (\ref{Qdependent}) satisfying the conditions (\ref{basic1c}) and (\ref{basic2c}).

We now come to solve its effective difference number.

{\bf Step 1:} First, we make the change of $\alpha=\theta+{1\over6}$ to transform (\ref{simpleexample1}) to
\begin{equation}
\label{88a}
\hat{\alpha}_{1}=-\alpha,\ \hat{\alpha}_{2}=-2\alpha+{2\over3},\ \hat{\alpha}_{3}=3\alpha.
\end{equation}
By Lemma 4.2, (\ref{88a}) is equivalent to
\begin{equation}
\label{88b}
\hat{\alpha}_{1}=-\alpha,\ \hat{\alpha}_{2}=-2\alpha+{2\over3},\ \hat{\alpha}_{31}=\alpha,\
\hat{\alpha}_{32}=\alpha+{1\over3},\ \hat{\alpha}_{33}=\alpha+{2\over3}.
\end{equation}

Secondly, we make the change of $\beta=\alpha-{1\over3}$ to transform (\ref{88b}) to
\begin{equation}
\label{88c}
\hat{\beta}_{1}=-\beta+{2\over3},\ \hat{\beta}_{2}=-2\beta,\ \hat{\beta}_{31}=\beta+{1\over3},\
\hat{\beta}_{32}=\beta+{2\over3},\ \hat{\beta}_{33}=\beta.
\end{equation}
Again by Lemma 4.2, (\ref{88c}) is equivalent to
\begin{equation}
\label{88d}
\hat{\beta}_{1}=-\beta+{2\over3},\ \hat{\beta}_{21}=-\beta,\ \hat{\beta}_{22}=-\beta+{1\over2},\
\hat{\beta}_{31}=\beta+{1\over3},\ \hat{\beta}_{32}=\beta+{2\over3},\ \hat{\beta}_{33}=\beta.
\end{equation}

{\bf Step 2:} By Lemma 4.3, we can cut off the following superfluous pairs in (\ref{88d}):
$$ \hat{\beta}_{1}=-\beta+{2\over3}\ \&\ \hat{\beta}_{31}=\beta+{1\over3};\quad \text{and}\quad  \hat{\beta}_{21}=-\beta\
\&\ \hat{\beta}_{33}=\beta. $$
That is,  (\ref{88d}) is equivalent to
\begin{equation}
\label{88e}
\hat{\beta}_{22}=-\beta+{1\over2},\ \ \hat{\beta}_{32}=\beta+{2\over3}.
\end{equation}
Finally, we make the change of $\gamma=\beta+{2\over3}$ to transform (\ref{88e})  to
\begin{equation}
\label{88f}
\hat{\gamma}_{22}=-\gamma+{1\over6},\ \ \hat{\gamma}_{32}=\gamma.
\end{equation}
It is obvious that the effective difference number of (\ref{88f}) is $1$ and so the system
(\ref{simpleexample1}) does. \hfill$\Box$


\section{Proof of Theorem \ref{mainresult}}

In this section, we prove our main Theorem \ref{mainresult}. Firstly
we give a proof of Theorem \ref{mainresult}  for $(\R P^{2n+1},F)$ which is involved in the irrational system
$\{\hat{\theta}_{1}$, $\hat{\theta}_{2}$, \dots, $\hat{\theta}_{k}\}$ with $2\leq k\leq 2n$ satisfying (\ref{mean}).
For sake of readability, we divide it into two cases according to whether $\rank(\hat{\theta}_{1},\hat{\theta}_{2},\dots\hat{\theta}_{k})=1$ or not.
We will give in details the proof for the first case. Based on the well known Kronecker's approximation theorem in Number theory,
the second one can be then proved quite similarly and so we only sketch it. While for $(\R P^{2n},F)$, the proof is similar and  will be explained at the end of this section.

{\bf Proof of Theorem \ref{mainresult} for $(\R P^{2n+1},F)$:} We carry out the proof into two cases.

{\bf Case 1:} $r=\rank(\hat{\theta}_{1},\hat{\theta}_{2},\dots\hat{\theta}_{k})=1.$

As we have mentioned in Section 4, the irrational system (\ref{mean}) with $r=1$ can be seen as a special case
of (\ref{Qdependent}) satisfying (\ref{basic1c}) and (\ref{basic2c}).

Since any $\eta$-action with $\eta\in\Q$  to (\ref{Qdependent}), if necessary, does no substantive effect on our following arguments,
by Theorem \ref{crucialtheorem} and Remark \ref{transformation} we can assume without loss of generality that
$$|k_{0}^{+}-k_{0}^{-}|\geq1\ \text{and}\ \mathcal{K}_{1}=\{1,2,\dots, k_{1}\},$$
with $k_{1}\geq1$ due to (\ref{basic2c}), and denote by $\xi_{j}={r_{j}\over q_{j}}$ for $1\leq j\leq k_{1}.$

 Let $\bar{q}=q_{1}q_{2}\cdots q_{k_1}$ and $m_{l}=2(n+1)\bar{q}l+1$ with $l\in\N$. Then, by (\ref{Qdependent}) we have
\begin{equation}
\label{sum3}
\begin{aligned}
\sum_{j=2}^{k}\left\{m_{l}\hat{\theta}_{j}\right\}&=\sum_{j=2}^{k_1}\left\{m_{l}\hat{\theta}_{j}\right\}+\sum_{j=k_{1}+1}^{k}\left\{m_{l}\hat{\theta}_{j}\right\}\\
&=\sum_{j=2}^{k_1}\left\{p_{j}\{m_{l}\theta\}+\xi_{j}\right\}+\sum_{j=k_{1}+1}^{k}\left\{p_{j}\{m_{l}\theta\}\right\}.\\
\end{aligned}
\end{equation}
Due to $\theta\in\mathbb{Q}^{c}$,  the set $\{\{m_{l}\theta\}\mid l\in\N\}$ is dense in $[0,1]$. For every $L\in\Z$, we introduce the  auxiliary function
\begin{equation}
\label{auxifunction}
f_{L}(x)=\sum_{j=2}^{k_1}\left\{\left\{p_{j}x+\xi_{j}\right\}+2L\hat{\theta}_{j}\right\}+\sum_{j=k_{1}+1}^{k}\left\{\left\{p_{j}x\right\}+2L\hat{\theta}_{j}\right\},\ \forall x\in [0,1],
\end{equation}
and denote for simplicity by $f=f_{0}$.

Let $a$ and $b$ in $(0,1)$ be two real numbers  sufficiently close to $0$ and $1$ respectively.
Then,
\begin{equation}
\label{fa}
\begin{aligned}
f(a)=\sum_{j=2}^{k_1}\left\{p_{j}a+\xi_{j}\right\}+\sum_{j=k_{1}+1}^{k}\left\{p_{j}a\right\}&=\sum_{j=2}^{k_1}\left(p_{j}a+\xi_{j}\right)+\sum_{j\in\mathcal{K}_{0}^{+}}p_{j}a+\sum_{j\in\mathcal{K}_{0}^{-}}\left(1+p_{j}a\right)\\
&=k_{0}^{-}+\sum_{j=2}^{k}p_{j}a+\sum_{j=2}^{k_1}\xi_{j},
\end{aligned}
\end{equation}
and by similar computation,
\begin{equation}
\label{fb}
\begin{aligned}
f(b)&=k_{0}^{+}+\sum_{j=2}^{k}p_{j}(b-1)+\sum_{j=2}^{k_1}\xi_{j}.
\end{aligned}
\end{equation}
It follows by (\ref{fa}) and (\ref{fb}) that
\begin{equation}
\label{distance}
|f(b)-f(a)|=|k_{0}^{+}-k_{0}^{-}+\sum_{j=2}^{k}p_{j}(b-1-a)|=|k_{0}^{+}-k_{0}^{-}+p_{1}(-b+1+a)|,
\end{equation}
where the second identity we have used $\sum_{j=1}^{k}p_{j}=0.$

\begin{lemma}
 \label{lemma1616a}
 If $a$ and $b$ in $(0,1)$ are sufficiently close to $0$ and $1$ respectively, then

{\rm (i)} $f(a)$ and $f(b)$ lie in different intervals of  (\ref{shortdengjia}) with $L=0$,

{\rm (ii)} $f_{L}(a)$ and $f_{L}(b)$ lie in the same interval of (\ref{shortdengjia}) for any $1\leq |L|\leq\bar{N}$ with $\bar{N}\in\N$ prior fixed, including $f_{L}(0).$
\end{lemma}
\Proof (i) By (\ref{distance}) and the assumption, $|f(b)-f(a)|\approx|k_{0}^{+}-k_{0}^{-}|.$
 Here and later,  we use $X\approx Y$ as usual to mean that $X$ is sufficiently close to $Y$ in the context of writing.
 Since the length of each interval in (\ref{shortdengjia}) with $L=0$ is less than or equal to $1$, so $f(a)$ and $f(b)$
 must lie in different ones, provided that  $|k_{0}^{+}-k_{0}^{-}|\geq2$.

 If $|k_{0}^{+}-k_{0}^{-}|=1$,  then $|f(b)-f(a)|\approx1$. For the case of $k=2$, since the length of each interval of  (\ref{shortdengjia}) with $L=0$ is less than $1$, (i) follows immediately. The rest case is $k\geq3$, which still contains three subcases.

$1^{\circ}$ If $k_{1}\geq2$, we get by (\ref{fa}) that
$$\{f(a)\}\approx\left\{k_{0}^{-}+\sum_{j=2}^{k_1}\xi_{j}\right\}=\left\{-\xi_{1}+\sum_{j=1}^{k_1}\xi_{j}\right\}
=\left\{\left\{\sum_{j=1}^{k_1}\xi_{j}\right\}-\xi_{1}\right\}=\left\{ \{Q_{0}\}-\xi_{1}\right\}.$$ Notice that the dividing points of the intervals in (\ref{shortdengjia}) with $L=0$ are $$0,\ \{Q_{0}\},\ 1+\{Q_{0}\},\ 2+\{Q_{0}\},\ \dots,\ k-2+\{Q_{0}\},\ k-1.$$ Therefore,
$\{k_{0}^{-}+\sum_{j=2}^{k_1}\xi_{j}\}=\left\{ \{Q_{0}\}-\xi_{1}\right\}$ must be an interior point of these intervals, so does $f(a)$. It then  yields that $f(a)$ and $f(b)$ must lie in two different intervals.

$2^{\circ}$ If $k_{1}=1$ and $k_{0}^{-}\geq1$, then  $f(a)\approx k_{0}^{-}$ is also an interior point and (i) follows.

$3^{\circ}$ If $k_{1}=1$ and $k_{0}^{-}=0$, then $f(a)=\sum_{j=2}^{k}p_{j}a$ lies in the first interval whose length is $\{Q_{0}\}<1$ and so $f(b)$ must lie in another one.

(ii) It can be checked directly that $\lim_{a\to0}f_{L}(a)=\lim_{b\to1}f_{L}(1)=f_{L}(0)\in\mathbb{Q}^{c},$ since $\xi_{j}\in\Q$ for $1\leq j\leq k$ and $\sum_{j=2}^{k}2L\hat{\theta}_{j}\in\mathbb{Q}^{c}$. But the dividing points of these intervals in (\ref{shortdengjia}) with $1\leq|L|\leq\bar{N}$ are finitely many rational numbers, so $f_{L}(0)$ is an interior point of these intervals and (ii) follows. \hfill$\Box$

Notice that $f$ is almost continuous on $(0,1)$. Without loss of generality, we assume $a$ and $b$ to be two points of continuity of $f$ and choose $l_{1}$, $l_{2}\in\N$ with $l_{2}-l_{1}$ sufficiently large such that $\{m_{l_1}\theta\}\approx a$ and $\{m_{l_2}\theta\}\approx b.$ Then by (\ref{sum3}), (\ref{auxifunction}) and (i) of Lemma \ref{lemma1616a}, we get
$\sum_{j=2}^{k}\left\{m_{l_1}\hat{\theta}_{j}\right\}$ and $\sum_{j=2}^{k}\left\{m_{l_2}\hat{\theta}_{j}\right\}$ lie in different intervals of (\ref{shortdengjia}) with $L=0$.  Suppose that $$\sum_{j=2}^{k}\left\{m_{l_1}\hat{\theta}_{j}\right\}\in I_{i'}\ \text{and}\ \sum_{j=2}^{k}\left\{m_{l_2}\hat{\theta}_{j}\right\}\in I_{i''},$$
with $\{i',i''\}\subseteq\{0,1,2,\dots,k-1\}$ and $i'\neq i''$. By (\ref{shortdengjia}) we have
$i(c^{m_{l_{1}}})=2n\bar{q}l_{1}+2[Q_{0}]-2i'$ and
\begin{equation}
\label{73a}
i(c^{m_{l_2}})=2n\bar{q}l_{2}+2[Q_{0}]-2i''.
\end{equation}
Since $2n\mid(2n\bar{q}l_{1}+2[Q_{0}]-2i'')$ if and only if $2n\mid(2n\bar{q}l_{2}+2[Q_{0}]-2i'')$, we get by (\ref{b.1}) that
 $$\beta_{2n\bar{q}l_{1}+2[Q_{0}]-2i''}=\beta_{2n\bar{q}l_{2}+2[Q_{0}]-2i''}\equiv\beta.$$
 Take $\bar{N}>4(n+1)$ in (ii) of Lemma \ref{lemma1616a} and observe $|2[Q_{0}]-2i''|\leq k\leq2n.$ By Lemma \ref{morsebetti} and Lemma \ref{bound}, there exist $L_{i}\in\Z$ with $1\leq |L_{i}|\leq \bar{N}$ and  $1\leq i\leq\beta$ such that
$$i(c^{m_{l_{1}}+2L_{i}})=2n\bar{q}l_{1}+2[Q_{0}]-2i''.$$
 Since $\sum_{j=2}^{k}\left\{(m_{l_{2}}+2L_{i})\hat{\theta}_{j}\right\}$ and
$\sum_{j=2}^{k}\left\{(m_{l_{1}}+2L_{i})\hat{\theta}_{j}\right\}$ are in the same interval of (\ref{shortdengjia}) with $1\leq |L_{i}|\leq\bar{N}$ by (ii) of Lemma \ref{lemma1616a}, we get again by (\ref{shortdengjia}) that
\begin{equation}
\label{73b}
i(c^{m_{l_{2}}+2L_{i}})=2n\bar{q}l_{2}+2[Q_{0}]-2i'',\ \forall 1\leq i\leq \beta.
\end{equation}
By (\ref{73a}) and (\ref{73b}), it yields $\beta\equiv\beta_{2n\bar{q}l_{2}+2[Q_{0}]-2i''}=\beta+1$ which is obviously absurd.

{\bf Case 2:} $r=\rank(\hat{\theta}_{1},\hat{\theta}_{2},\dots\hat{\theta}_{k})\geq2.$

By Lemma \ref{L3.1},  there are $p_{jl}\in\Z$, $\theta_{k_l}\in\Q^{c}$ and $\xi_{j}\in\Q$ with $1\leq l\leq r$ and $1\leq j\leq k$ such that
\begin{equation}
\label{system82a}
\hat{\theta}_{j}=\sum_{l=1}^{r}p_{jl}\theta_{k_l}+\xi_{j},\ \forall 1\leq j\leq k.
\end{equation}
Moreover,  $\theta_{k_1}$, $\theta_{k_2}$ \dots, $\theta_{k_r}$ are linearly independent over $\Q$. Due to
(\ref{mean}), it follows
\begin{equation}
\label{20150923a}
\sum_{j=1}^{k}p_{jl}=0,\ \forall 1\leq l\leq r.
\end{equation}.

Our basic idea for proving Case 2 is to construct an irrational system  with rank $1$  associated to (\ref{system82a}), which plays the essential role in our sequel arguments  due to  the following result.\\
{\bf Kronecker's approximation theorem} (cf. Theorem 7.10 in \cite{Apostol}): {\it If $\theta_{1}$, $\theta_{2}$, \dots, $\theta_{r}$ are linearly independent over $\Q$, then the set $\{(m\theta_{1},m\theta_{2},\dots,m\theta_{r})\mid m\in N\}$ is dense in
$$[0,1]^{r}=\underbrace{[0,1]\times[0,1]\times\cdots\times[0,1]}_{r}.$$}
\begin{lemma}
\label{lemma82a}
There are $s_{2},\ s_{3}, \dots,\ s_{r}\in\Z$ such that
\begin{equation}
\label{82a}
p_{j1}+\sum_{l=2}^{r}s_{l}p_{jl}\in\Z\backslash\{0\},\ \forall 1\leq j\leq k,
\end{equation}
\end{lemma}
\Proof  Let $J_{0}=\{1\leq j\leq k\mid p_{j1}=0\}.$ If $J_{0}=\emptyset$, we need only take $s_{2}=s_{3}=\cdots=s_{r}=0.$
 If $J_{0}\neq\emptyset$,  we claim that $(p_{j2},p_{j3},\dots,p_{jr})\neq(0,0,\dots,0)$ for each $j\in J_{0}$. Otherwise,  then (\ref{system82a}) yields that
$\hat{\theta}_{j}=\xi_{j}\in\Q,$
which contradicts to $\hat{\theta}_{j}\in\Q^{c}$. So the set
\begin{equation*}
\label{82b}
X_{j}\equiv\left\{(x_{2},x_{3},\dots,x_{r})\mid p_{j2}x_{2}+p_{j3}x_{3}+\dots+p_{jr}x_{r}=0\right\},
\end{equation*}
 is a subspace of dimension $r-2$ in $\R^{r-1}$ which yields that $X=\cup_{j\in J_{0}}X_{j}$ is a proper subset of $\R^{r-1}$. Take arbitrarily out an integral point $(\bar{s}_{2},\bar{s}_{3},\dots,\bar{s}_{r})\in\R^{r-1}\backslash X$.
Then for every $\bar{N}\in\N$ we have
\begin{equation}
\label{20150922a}
|p_{j1}+\sum_{l=2}^{r}\bar{N}\bar{s}_{l}p_{jl}|
=\left\{\begin{array}{ll}
\bar{N}|\sum_{l=2}^{r}\bar{s}_{l}p_{jl}|\neq0,& \text{if}\ j\in J_{0},\\
|p_{j1}|\neq0,& \text{if}\ j\notin J_{0}\ \text{and}\ \sum_{l=2}^{r}\bar{s}_{l}p_{jl}=0,\\
|p_{j1}+\bar{N}\sum_{l=2}^{r}\bar{s}_{l}p_{jl}|,& \text{if}\ j\notin J_{0}\ \text{and}\ \sum_{l=2}^{r}\bar{s}_{l}p_{jl}\neq0.
\end{array}
\right.
\end{equation}
For the third case in the righthand side of (\ref{20150922a}), we can take $\bar{N}\in\N$ sufficiently large so that $|p_{j1}+\bar{N}\sum_{l=2}^{r}\bar{s}_{l}p_{jl}|\neq0$ for all these $j$'s therein. Finally let $s_{l}=\bar{N}\bar{s}_{l}$  and (\ref{82a}) follows. \hfill$\Box$

By Lemma \ref{lemma82a}, we can make the change of variables
$\td{\theta}_{k_1}=\theta_{k_1}\ \text{and}\ \td{\theta}_{k_l}=\theta_{k_l}-s_{l}\theta_{k_1}\ \text{for}\ 2\leq l\leq r.$
Then the system (\ref{system82a}) is transformed to
\begin{equation}
\label{system82b}
\hat{\theta}_{j}=\sum_{l=1}^{r}\td{p}_{jl}\td{\theta}_{k_l}+\xi_{j},\ \forall 1\leq j\leq k,
\end{equation}
with
$\td{p}_{j1}=p_{j1}+\sum_{l=2}^{r}s_{l}p_{jl}\in\Z\backslash\{0\},$ and by (\ref{20150923a}) we have
\begin{equation*}
\label{basic2ccc}
\sum_{j=1}^{k}\td{p}_{j1}=\sum_{j=1}^{k}p_{j1}+\sum_{j=1}^{k}\sum_{l=2}^{r}s_{l}p_{jl}=0+\sum_{l=2}^{r}s_{l}\left(\sum_{j=1}^{k}p_{jl}\right)=0.
\end{equation*}
Since  $\theta_{k_1}$, $\theta_{k_2}$, \dots, $\theta_{k_r}$ are linearly independent over $\Q$, so do  $\td{\theta}_{k_1}$, $\td{\theta}_{k_2}$, \dots, $\td{\theta}_{k_r}.$

Consider the following irrational system  with rank $1$ associated to (\ref{system82b})
\begin{equation}
\label{system82c}
\hat{\alpha}_{j}=\td{p}_{j1}\td{\theta}_{k_1}+\xi_{j},\ \forall 1\leq j\leq k.
\end{equation}
By Theorem \ref{crucialtheorem}, without loss of generality we can assume for (\ref{system82c}) that $|\td{k}_{0}^{+}-\td{k}_{0}^{-}|\geq1$ and
$\td{\mathcal{K}}_{1}=\{1,2,\dots, \td{k}_{1}\}$, and denote by $\xi_{j}={r_{j}\over q_{j}}$ for $1\leq j\leq \td{k}_{1}.$

 Let $\td{q}=q_{1}q_{2}\cdots q_{\td{k}_{1}}$ and $\td{m}_{l}=2(n+1)\td{q}l+1$ for $l\in\N$. Then, we get by (\ref{system82b}) that
\begin{equation}
\label{sum33}
\begin{aligned}
\sum_{j=2}^{k}\left\{\td{m}_{l}\hat{\theta}_{j}\right\}
&=\sum_{j=2}^{\td{k}_1}\left\{\td{m}_{l}\hat{\theta}_{j}\right\}+\sum_{j=\td{k}_{1}+1}^{k}\left\{\td{m}_{l}\hat{\theta}_{j}\right\}\\
&=\sum_{j=2}^{\td{k}_1}\left\{\sum_{l=1}^{r}\td{p}_{jl}\{\td{m}_{l}\td{\theta}_{k_l}\}+\xi_{j}\right\}+\sum_{j=\td{k}_{1}+1}^{k}\left\{\sum_{l=1}^{r}\td{p}_{jl}\{\td{m}_{l}\td{\theta}_{k_l}\}\right\}\\
\end{aligned}
\end{equation}

By Kronecker's approximation theorem,  the set $\{(\{\td{m}_{l}\td{\theta}_{k_1}\},\{\td{m}_{l}\td{\theta}_{k_2}\},\dots,\{\td{m}_{l}\td{\theta}_{k_r}\})\mid l\in\N\}$ is dense in $[0,1]^{r}$. For every $L\in\Z$, we can introduce the  auxiliary multi-variable function on $[0,1]^{r}$,
$$g_{L}(x_{1},x_{2},\dots,x_{r})=\sum_{j=2}^{\td{k}_1}\left\{\sum_{l=1}^{r}\td{p}_{jl}x_{l}+\xi_{j}+2L\hat{\theta}_{j}\right\}
+\sum_{j=\td{k}_{1}+1}^{k}\left\{\sum_{l=1}^{r}\td{p}_{jl}x_{l}+2L\hat{\theta}_{j}\right\},$$
and denote for simplicity by $g=g_{0}.$ Similarly as before, we have
 \begin{lemma}
\label{lemma1616b}
  If $(a_{1},a_{2},\dots, a_{r})$ and $(b_{1},b_{2},\dots, b_{r})$ in $(0,1)^{r}$ are sufficiently close to $(0,0,0,\dots,0)$ and $(1,0,0,\dots,0)$ respectively by a suitable means, then

{\rm (i)} $g(a_{1},a_{2},\dots, a_{r})$ and $g(b_{1},b_{2},\dots, b_{r})$ lie in different intervals of  (\ref{shortdengjia}) with $L=0$.

 {\rm (ii)} $g_{L}(a_{1},a_{2},\dots, a_{r})$ and $g_{L}(b_{1},b_{2},\dots, b_{r})$
lie in the same interval of (\ref{shortdengjia}) for any $1\leq |L|\leq\bar{N}$ with $\bar{N}\in\N$ prior fixed, including $g_{L}(0,0,\dots, 0).$
\end{lemma}
\Proof
 (i) Since $a_1$, $a_2$, \dots, $a_r$ (resp. $b_1$, $b_2$, \dots, $b_r$) are independent, we can select them by such a way that the decimal functions in $g(a_{1},a_{2},\dots, a_{r})$ and $g(b_{1},b_{2},\dots, b_{r})$ are mainly determined by $a_1$ and $b_1$ respectively.  For instance, this can be realized by requiring $a_l$ (resp. $b_l$) with $2\leq l\leq r$ to be much smaller than $a_1$ (resp. $1-b_1$).  The rest proof is then similar as that in Lemma \ref{lemma1616a}-(i), with $g$ in stead of $f$ therein.

 (ii) follows the same line as  Lemma \ref{lemma1616a}-(ii) and do not need such a choice as above.

Due to Lemma \ref{lemma1616b}, the rest proof is then almost word by word as that in Case 1 and so we omit the tedious details. \hfill$\Box$

\begin{remark}
As for $\R P^{3}$, we can give a more direct and easier proof. Indeed, we can make a reduction by (\ref{mean}) (with  $n=1$ and $k=2$) so that only one irrational number is rest.
The uniformly distribution mod one in Number theory then enables the authors in \cite{DLX2015} to find some $l\in\N$
such that the Betti number $\bar{\beta}_{2l}=1$ which contradicts to the topological structure of the non-contractible
loop space on $\R P^{3}$ obtained in \cite{XL2015}. However when one tries to use such a means  to deal with higher
dimensional $\R P^{2n+1}$, more irrational numbers are rest to be controlled simultaneously for larger $k.$ What is even worse, those irrational
numbers may be linearly dependent over $\Q$. These facts make the arguments in section 3.3 of \cite{DLX2015} difficult to continue,
even for $\R P^5$.

 To overcome these difficulties, we discover a general character of the irrational systems (\ref{Qdependent}) satisfying the conditions (\ref{basic1c}) and (\ref{basic2c}), which are closely associated to our problem. That is, the effective difference number of each of such irrational systems is larger than or equal to $1$ (cf. Theorem \ref{crucialtheorem}). Based on it and the Kronecker's approximation theorem, we can get the desired contradiction dynamically (quite different from the static way in [14]), provided that there is only one non-contractible closed geodesics.
\end{remark}

{\bf Proof of Theorem \ref{mainresult} for $(\R P^{2n},F)$:}

For the case of even $n$, it shares the same essential properties with the odd case except for
some quantitative differences, such as the resonance identity, the precise index iteration formulae and the irrational systems.
Hence we only sketch its proof to avoid this paper being too long and tedious.

We now give some explanations to the proof of Theorem 1.2 for the case of $\R P^{2n}$.

Note that in the proof of Theorem 1.2 for the case of $\R P^{2n+1}$,  only (\ref{b.1}), (\ref{beqm})-(\ref{shortdengjia}), Lemmas \ref{morsebetti}-\ref{bound}
and Theorem \ref{crucialtheorem} are used. As for the case of $\R P^{2n}$, Lemma \ref{morsebetti} and
Theorem \ref{crucialtheorem} still hold, thus by using (\ref{b.2}) instead of (\ref{b.1}), Lemma \ref{lemma3.4} instead of Lemma \ref{lemma3.2},
Lemma \ref{bound2} instead of Lemma \ref{bound}, and noticing  Remark 3.2,
we can go through the proof of Theorem 1.2 for the case of $\R P^{2n}$ word by word as that of the case of $\R P^{2n+1}$.
We complete the proof of Theorem 1.2.\hfill$\Box$

\smallskip{\noindent\bf Acknowledgements.}
We would like to thank sincerely the referee for his/her careful reading, valuable comments
on this paper, and for his/her deep insight on the main ideas of this paper. And also we would like to sincerely thank our
advisor, Professor Yiming Long, for introducing us to the theory of
closed geodesics and for his valuable advices and comments on this paper.

\end{document}